\documentstyle[amssymb]{amsart}

\def\fun#1#2{\lower3.6pt\vbox{\baselineskip0pt\lineskip.9pt
  \ialign{$\mathsurround=0pt#1\hfil##\hfil$\crcr#2\crcr\sim\crcr}}}
\newskip\humongous \humongous=0pt plus 1000pt minus 1000pt

\newif\ifdtup

\def\oldreffmt#1{\rlap{[#1]} \hbox to 2\parindent{}}

\def\figfmt#1{\rlap{Figure {#1}} \hbox to 1in{}}

\def\beq{\begin{equation}}
\def\eeq{\end{equation}}

\def\bq{\begin{quote}}
\def\eq{\end{quote}}
\def\lqqp{\lq \lq \ }
\def\lqq{\lq \lq }
\relax

\newtheorem{THEOREM}{Theorem}[section]
\newtheorem{LEMMA}{Lemma}[section]

\theoremstyle{definition}
\newtheorem{definition}{Definition}[section]

\newcommand{\be}{\begin{equation}}
\newcommand{\ee}{\end{equation}}
\addtocounter{equation}{-1}
\begin{document}
\rightline{In {\it Lie Algebras,  Cohomologies  and  New  Findings  in  Quantum Mechanics}}
\rightline{AMS `Contemporary Mathematics' series, N. Kamran and P. Olver (Eds.),}
\rightline{vol. 160, pp. 263-310 (1994) }
\vskip .3truecm
\title[Lie-algebras and linear operators \dots]
{Lie-algebras and linear operators with invariant subspaces }
\author{ Alexander Turbiner}
\address
{Institute for Theoretical and Experimental Physics,
Moscow 117259, Russia}
\email{TURBINER@@VXCERN.CERN.CH or TURBINER@@CERNVM}
\curraddr{I.H.E.S., Bures-sur-Yvette, F-91440, France}
\thanks{Supported in Part by the Swiss National Foundation.}
\thanks{This paper is in final form
and no version will be submitted for publication elsewhere.}
\thanks{The manuscript is prepared using {\em AMS-LATEX}.}

\date{}
\maketitle
\begin{abstract}
A general classification of linear differential
 and finite-difference operators possessing a finite-dimensional
invariant subspace with a polynomial basis (the generalized Bochner
problem) is given. The main result is that
any operator with the above property must have a representation as a
polynomial element of the universal enveloping algebra of
some algebra of differential (difference) operators in
finite-dimensional representation plus an operator annihilating
the finite-dimensional invariant subspace. In low dimensions a
classification is given by algebras
$sl_2({\bold R})$ (for differential operators in ${\bold R}$)
and $sl_2({\bold R})_q$ (for finite-difference operators in  ${\bold R}$),
 $osp(2,2)$  (operators in one real and one Grassmann variable, or
equivalently, $2 \times 2$ matrix operators in ${\bold R}$),
$sl_3({\bold R})$,  $sl_2({\bold R}) \oplus sl_2({\bold R})$
and  $gl_2 ({\bold R}) \ltimes {\bold R}^{r+1}\ , r$ a
natural number (operators in ${\bold R^2}$). A
 classification of linear operators possessing infinitely many
finite-dimensional invariant subspaces with a basis in polynomials
is presented. A connection to the recently-discovered
quasi-exactly-solvable spectral problems is discussed.
\end{abstract}

 S. Bochner (1929) asked about a classification of differential equations
\begin{equation}
 T \varphi  \  = \ \epsilon \varphi
\end{equation}
where $T$ is a linear differential operator of $k$th order in one real
variable $x \in {\bold R}$
and $\epsilon$ is the spectral parameter, having an infinite sequence of
{\it orthogonal} polynomial solutions \cite{b} (see also \cite{little}).
\begin{definition} Let us give the name of the {\it generalized Bochner
 problem} to the problem of classification of linear differential
(difference) operators, for which the eigenvalue problem
 (0) has a certain number of eigenfunctions in the form of a
finite-order polynomial in some variables.
\end{definition}

Following this definition the original Bochner problem is simply a
particular case.
In \cite{t1} \ a general method has been formulated for generating
 eigenvalue problems for linear differential
 operators, linear matrix differential operators and
 linear finite-difference operators in one and several variables
possessing polynomial solutions. The
method was based on considering the eigenvalue problem for the
representation of a polynomial element
of the universal enveloping algebra of the Lie algebra in a
finite-dimensional, \lq projectivized'
representation of this Lie algebra \cite{t1}. Below it is shown that
 this method provides both necessary and sufficient conditions for
the existence of polynomial solutions of linear differential equations
and a certain class of finite-difference  equations.

The generalized Bochner problem can be subdivided into two parts: (i)
a classification of linear operators possessing an invariant subspace
(subspaces) with a basis in polynomials of finite degree and
(ii) a description of the
conditions under which such operators are symmetrical.
This paper  will be devoted to a solution of the first problem; as for
the second one, only the first step has been done (see below).
The plan of the paper is the following: Section 1 is devoted to the
case of differential operators in one real variable; in Section 2
finite-difference operators in ${\bold R}$ are treated; operators
in one real and one Grassmann variables are considered in Section 3,
while the case of operators in  ${\bold R^2}$ is given in Section 4.
The general situation is described in the Conclusion.
\section{ Ordinary differential equations}
Consider the space of all polynomials of order $n$
\begin{equation}
\label{e1}
{\cal P}_n \ = \ \langle 1, x, x^2, \dots , x^n \rangle ,
\end{equation}
where $n$ is a non-negative integer and $x \in {\bold R}$.
\begin{definition} Let us name a linear differential operator of the $k$th
order, $T_k$ {\it quasi-exactly-solvable}, if it preserves the space
${\cal P}_n$. Correspondingly, the operator $E_k$, which preserves
the infinite flag $ {\cal P}_0 \subset  {\cal P}_1 \subset {\cal P}_2
\subset \dots \subset {\cal P}_n \subset \dots $ of spaces of all
polynomials, is named {\it exactly-solvable}.
\end{definition}
\begin{LEMMA} {\it (i) Suppose $n > (k-1)$.  Any quasi-exactly-solvable
operator $T_k$, can be
represented by a $k$-th degree polynomial of the operators}
\label{e2}
\[ J^+ = x^2 \partial_x - n x,\  \]
\begin{equation}
 J^0 = x \partial_x - {n \over 2} \  ,
\end{equation}
\[ J^- = \partial_x \ ,  \]
\noindent
{\it (the operators (2) obey the $sl_2({\bold R})$ commutation relations}
\footnote{The representation (2) is one of the \lq projectivized'
representations (see \cite{t1}). This realization of $sl_2({\bold R})$
has been derived at the first time by Sophus Lie.}
{\it). If $n \leq (k-1)$, the part of the quasi-exactly-solvable operator
$T_k$ containing
derivatives up to order $n$ can be represented by an $n$th
degree polynomial in the generators (2).

(ii) Conversely, any polynomial in (2) is quasi-exactly solvable.

(iii) Among quasi-exactly-solvable operators
there exist exactly-solvable operators $E_k \subset T_k$.}
\end{LEMMA}

{\it Comment 1.} If we define the universal enveloping algebra $U_g$ of a
Lie algebra $g$ as the algebra of all ordered polynomials in generators, then
$T_k$ at $k < n+1$ is simply an element
 of the universal enveloping algebra $U_{sl_2({\bold R})}$ of the algebra
$sl_2({\bold R})$ taken in representation (2). If $k \geq n+1$, then $T_k$ is
represented as an element of $U_{sl_2({\bold R})}$ plus $B {d^{n+1}
\over dx^{n+1}}$, where $B$ is any linear differential operator of order
not higher than $(k-n-1)$.
\begin{pf}
The essential part of the proof is based on the Burnside theorem
(see, e.g., \cite{w}\footnote{Burnside theorem is a particular case
of more general Jacobson theorem (see \cite{l}, Chapter XVII.3).
I am grateful to V. Kac for this comment.}):
\begin{quote}
Let $A_1,\dots,A_k$ be linear operators in a linear space $E$ over real
numbers, $dimE <\infty$. Let us assume that there is no linear space
$L$ over real numbers,\ $0<dimL <dimE$ , such that $A_i: L \mapsto L$
for all $i=1,2,\ldots,k$.
Then any linear operator acting in $E$ can be represented as a
polynomial in $A_1,\dots,A_k$.
\end{quote}
The operators $J^{\pm,0}$ act in ${\cal P}_n$ irreducibly and therefore
the theorem can be applied. Also there exist operators $B$ having
${\cal P}_n$ as a kernel, $B: {\cal P}_n \mapsto 0$. Clearly, those operators
have a form $B(x,\partial_x) {d^{n+1} \over dx^{n+1}}$,
where $B(x,\partial_x)$ is any linear differential operator, and they make
no contribution in $T_k$ if $k < n+1$. It completes the proof of part (i).
Parts (ii) and (iii) are easy to prove based on part (i).\footnote{A detailed,
technical proof is presented in \cite{t2}.}
\end{pf}

Since $sl_2({\bold R})$ is a graded algebra, let us introduce the grading
of generators (2):
\begin{equation}
\label{e3}
deg (J^+) = +1 \ , \ deg (J^0) = 0 \ , \ deg (J^-) = -1 ,
\end{equation}
 hence
\begin{equation}
\label{e4}
deg [(J^+)^{n_+} (J^0)^{n_0}(J^-)^{n_-}] \  = \ n_+ - n_- .
\end{equation}
The grading allows us to classify the operators $T_k$ in the
Lie-algebraic sense.
\begin{LEMMA} {\it A quasi-exactly-solvable operator $T_k \subset
U_{sl_2({\bold R})}$ has no terms of positive grading, if and only if
it is an exactly-solvable operator.}
\end{LEMMA}
\begin{THEOREM} {\it Let $n$ be a non-negative integer. Take the eigenvalue
problem for a linear differential
operator of the $k$th order in one variable
\begin{equation}
\label{e5}
 T_k \varphi \ = \ \varepsilon \varphi \ ,
\end{equation}
where $T_k$ is symmetric. The problem (5) has $(n+1)$ linearly independent
eigenfunctions in the form of a polynomial in variable $x$ of order
not higher than $n$, if and only if $T_k$ is quasi-exactly-solvable.
The problem (5) has an infinite sequence of polynomial eigenfunctions,
if and only if the operator is exactly-solvable.}
\end{THEOREM}

\noindent
{\it Comment 2.} The \lqqp if "
 part of the first statement is obvious.
The \lqqp only if " part is a direct corollary of Lemma 1.1 .

This theorem gives a general classification of differential equations
\begin{equation}
\label{e6}
 \sum_{j=0}^{k} a_j (x) \varphi^{(j)} (x) \ = \ \varepsilon \varphi(x)
\end{equation}
having at least one polynomial solution in $x$.
The coefficient functions $a_j (x)$ must have the form
\begin{equation}
\label{e7}
a_j (x) \ = \ \sum_{i=0}^{k+j} a_{j,i} x^i
\end{equation}
The explicit expressions (7) for coefficient function in (6) are obtained by
the substitution (2) into a general, $k$th degree
 polynomial element of the universal
enveloping algebra $U_{sl_2({\bold R})}$.  Thus the coefficients $a_{j,i}$
can be expressed  through the coefficients of the
$k$th degree polynomial element of the universal
enveloping algebra $U_{sl_2({\bold R})}$. The number of free parameters of the
polynomial solutions is defined by the number of parameters
characterizing a general $k$th degree polynomial element of the universal
enveloping algebra $U_{sl_2({\bold R})}$. A rather straightforward calculation
leads to the following formula
\begin{equation}
\label{8}
 par (T_k) = (k+1)^2
\end{equation}
where we denote the number of free parameters of operator $T_k$ by the symbol
$par(T_k)$.
For the case of an infinite sequence of polynomial solutions expression
(7) simplifies to
\begin{equation}
\label{e9}
a_j (x) \ = \ \sum_{i=0}^{j} a_{j,i} x^i
\end{equation}
in agreement with the results of H.L. Krall's classification theorem
\cite{kr}(see also \cite{little}). In this case the number of free parameters
is equal to
\begin{equation}
\label{10}
 par (E_k) = {(k+1)(k+2) \over 2}
\end{equation}
In the present approach Krall's theorem is simply a description of
differential operators of $k$th order in one variable preserving
a finite flag  $ {\cal P}_0 \subset
 {\cal P}_1 \subset {\cal P}_2 \subset \dots \subset {\cal P}_k $ of spaces
of polynomials. One can easily show that the preservation of such a set of
polynomial spaces implies the preservation of an infinite flag of such spaces.

One may ask a more general question: which non-degenerate linear differential
operators have a finite-dimensional invariant sub-space of the form
\begin{equation}
\label{11}
\langle \alpha (z), \alpha (z) \beta (z), \dots , \alpha (z) \beta (z) ^n
\rangle \ ,
\end{equation}
where $\alpha (z),\beta (z)$ are some functions. Such operators
are easily obtained from the quasi-exactly-solvable operators
(see Lemma 1.1)
making the change of variable $x = \beta (z)$ and the \lqq gauge"
transformation $\tilde T = \alpha (z) T  \alpha (z)^{-1}$.
Since any one- or two-dimensional invariant sub-space can be presented in
the form (11), the important general statement takes place:
{\it there are no linear operators possessing one- or two-dimensional
invariant sub-space with an explicit basis other than given by
\linebreak Lemma 1.1} .
Therefore for any linear operator possessing one- or two-dimensional
invariant sub-space with an explicit basis the eigenvalue problem (5)
can be reduced to the form (16) or (15), respectively (see below).

Let us consider the second-order differential equation (6),
which can possess polynomial
solutions. From Theorem 1.1 it follows that the corresponding differential
operator must be quasi-exactly-solvable and can be represented as
\[ T_2 =  c_{++} J^+ J^+ + c_{+0} J^+ J^0 + c_{+-} J^+ J^- + c_{0-} J^0 J^-
+ c_{--} J^- J^- + \]
\begin{equation}
\label{e12}
 c_+ J^+ + c_0 J^0 + c_- J^- + c ,
\end{equation}
where $c_{\alpha \beta}, c_{\alpha}, c \in {\bf R}$.
The number of free parameters is $par (T_2) = 9$. Under the condition
$c_{++}  = c_{+0}  = c_+  =0$, the operator $T_2$ becomes
exactly-solvable (see Lemma 1.2) and the number of free parameters is
reduced to $par (E_2) = 6$.
\begin{LEMMA} {\it If the operator (12) is such that
\begin{equation}
\label{e13}
c_{++}=0 \quad and \quad c_{+} = ({n \over 2} - m)  c_{+0} \ , \ at \ some
\ m=0,1,2,\dots
\end{equation}
 then the operator $T_2$ preserves both ${\cal P}_n$ and
${\cal P}_m$. In this case the number of free parameters is
$par (T_2) = 7$.}
\end{LEMMA}

In fact, Lemma 1.3 means that $T_2 (J^{\alpha}(n),c_{\alpha \beta},
 c_{\alpha})$
can be rewritten as $T_2 (J^{\alpha}(m),\allowbreak c'_{\alpha
 \beta},c'_{\alpha})$.
As a consequence of Lemma 1.3 and Theorem 1.1, in general,
among polynomial solutions of (6) there
are polynomials of order $n$ and order $m$.

{\bf Remark.}  From the Lie-algebraic point of view Lemma 1.3 means
the existence of representations of second-degree polynomials in the
generators (2) possessing two invariant sub-spaces.
In general, if $n$ in (2) is a non-negative integer, then
 among representations of $k$th degree polynomials in the generators (2),
lying in the universal enveloping algebra, there exist representations
possessing $1,2,...,k$ invariant sub-spaces. Even starting from an
infinite-dimensional representation of the original algebra
($n$ in (2) is {\it not} a non-negative integer), one can construct
the elements of the universal enveloping algebra having finite-dimensional
representation (e.g., the parameter $n$ in (13) is non-integer,
however $T_2$ has the invariant sub-space of dimension $(m+1)$).
Also this property implies the existence of representations of
the polynomial elements of the universal enveloping algebra,
which can be obtained starting
from different representations of the original algebra.

Substituting (2) into (12) and then into (6), we obtain
\begin{equation}
\label{e14}
- P_{4}(x) \partial_x ^2 \varphi (x) \ +\ P_{3}(x) \partial_x  \varphi (x) \
+\ P_{2}(x) \varphi (x) \ =\ \varepsilon \varphi (x) ,
\end{equation}
where the $P_{j}(x)$ are polynomials of $j$th order with coefficients
related to  $ c_{\alpha \beta}, c_{\alpha}$ and $n$. In general,
problem (14) has $(n+1)$ polynomial solutions. If $n=1$, as a
consequence of Lemma 1.1, a  more
general spectral problem than (14) arises
\begin{equation}
\label{e15}
- F_{3}(x) \partial_x ^2 \varphi (x) \ +\ Q_{2}(x) \partial_x  \varphi (x) \
+\ Q_{1}(x) \varphi (x) \ =\ \varepsilon \varphi (x) ,
\end{equation}
where $F_3$ is an arbitrary real function of $x$ and $Q_j (x), j=1,2$ are
polynomials of order $j$, possessing only two polynomial solutions of the
form $(ax+b)$. For the case $n=0$ (one polynomial solution)
the spectral problem (6) becomes
\begin{equation}
\label{e16}
- F_{2}(x) \partial_x ^2 \varphi (x) \ +\ F_{1}(x) \partial_x  \varphi (x) \
+\ Q_0 \varphi (x) \ =\ \varepsilon \varphi (x) ,
\end{equation}
where $F_{2,1}(x)$ are arbitrary real functions and $Q_0$
is a real constant. After the transformation
\begin{equation}
\label{e17}
\Psi (z) \ =\ \varphi (x(z)) e ^ {-A(z)} ,
\end{equation}
where $x=x(z)$ is a change of the variable and $A(z)$
is a certain real function, one can reduce (14)--(16) to the
Sturm-Liouville-type problem
\begin{equation}
\label{e18}
(- \partial_z ^2 \ + \ V(z) ) \Psi (z) \ = \ \varepsilon \Psi (z),
\end{equation}
with the potential, which is equal to
\[ V(z) = (A')^2 - A'' + P_2 (x(z)) \ , \]
if
\[ A = \int ({P_3 \over P_4})dx - log z'\ ,
\ z = \int {dx \over \sqrt{P_4}}\ .\]
for the case of (14). If the functions (17), obtained after transformation,
belong to the ${\cal L}_2(\cal D)$-space\footnote{Depending on the
change of variable  $x = x(z)$, the space $\cal D$ can be whole
real line, semi-line and a finite interval.} ,
we reproduce the recently discovered quasi-exactly-solvable problems
\cite{t4}, where a finite number of eigenstates was found algebraically.
For example,
\begin{equation}
\label{e19}
 T_2 = -4 J^0 J^- + 4a J^+ + 4b J^0 - 2(n+1+2k) J^- + b(2n+1+2k)
\end{equation}
leads to the spectral problem (18) ($x=z^2$) with the potential
\begin{equation}
\label{e20}
 V(z) = a^2z^6 + 2abz^4 + [b^2 - (4n+3+2k)a]z^2 ,
\end{equation}
for which at $k=0 (k=1)$ the first $(n+1)$ eigenfunctions, even (odd)
in $x$, can be found algebraically. Of course, the number of those
`algebraized' eigenfunctions is nothing but the dimension of the
irreducible representation (1).

Taking different exactly-solvable operators $E_2$ for the eigenvalue
problem (6)
one can reproduce the equations having the Hermite, Laguerre, Legendre
and Jacobi polynomials as solutions \cite{t1}.\footnote{For instance,
setting the parameter $a=0$ in (19), the equation (15) converts to the
Hermite equation (after some substitution).} It is worth noting that
for general exactly-solvable operator $E_2$ the eigenvalues are
quadratic in number of eigenstate and can be presented as follows
\[
{\epsilon}_n =  \ c_{00} n^2 \ + \ c_{0} n \ + \ const
\]
(for details see  \cite{t1}).

Under special choices
of the general element $E_4 (E_6, E_8)$, one can reproduce all
known fourth-(sixth-, eighth-)order differential equations giving rise to
infinite sequences of orthogonal polynomials (see \cite{little} and
other papers in this volume).

  Recently, A. Gonz\'alez-Lop\'ez, N. Kamran and P. Olver \cite{olver}
gave the complete description of second-order polynomial elements of
$U_{sl_2({\bf R})}$ leading to the square-integrable eigenfunctions of the
Sturm-Liouville problem (18) after transformation (17). Consequently,
for second-order ordinary
differential equations (14) the combination of their result and Theorem 1.1
 gives a general solution of the Bochner problem  as well as the
more general problem of classification of  equations possessing a finite
number of orthogonal polynomial solutions.
\section{ Finite-difference equations in one variable}
Let us introduce the finite-difference analogue of the differential
operators (2) \cite{ot}
\label{e21}
\[ \tilde  J^+ = x^2 D - \{ n \} x \]
\begin{equation}
\tilde  J^0 = \  x D - \hat{n}
\end{equation}
\[ \tilde  J^- = \ D , \]
where $\hat n \equiv {\{n\}\{n+1\}\over \{2n+2\}}$ ,
$\{n\} = {{1 - q^n}\over {1 - q}}$
is the quantum symbol (or $q$-number), $q \in {\bold R}$ is a number
characterizing the deformation, $Df(x) = {{f(x) - f(qx)} \over
{(1 - q) x}}$ is a shift or a
finite-difference operator (or the so-called Jackson symbol (see \cite{e})).
The operators (21) after multiplication by some factors
\[ \tilde  j^0 = {q^{-n} \over q+1} {\{2n+2\} \over \{n+1\}} \tilde J^0 \]
\[ \tilde  j^{\pm} = q^{-n/2} \tilde  J^{\pm} \]
( see \cite{ot}) form a quantum $sl_2({\bold R})_q$ algebra with the following
commutation relations
\label{22}
\[ q \tilde  j^0\tilde  j^- \ - \ \tilde  j^-\tilde  j^0 \
= \ - \tilde  j^-  \]
\begin{equation}
 q^2 \tilde  j^+\tilde  j^- \ - \ \tilde  j^-\tilde  j^+ \
= \ - (q+1) \tilde  j^0
\end{equation}
\[ \tilde  j^0\tilde  j^+ \ - \ q\tilde  j^+\tilde  j^0 \ = \  \tilde  j^+  \]
(this algebra corresponds to the second Witten quantum deformation
of $sl_2$ in the classification of C. Zachos \cite{z}).
If $q \rightarrow 1$, the commutation relations (22) reduce
to the standard $sl_2({\bold R})$ ones. A remarkable property of generators
(21) is that, if $n$ is a non-negative integer, they form
the finite-dimensional representation. For $q$ others than root of unity
this representation is irreducible.

Similarly as for differential operators one can introduce
quasi-exactly-solvable $\tilde  T_k$ and exactly-solvable operators
 $\tilde  E_k$.
\begin{LEMMA} {\it (i) Suppose $n > (k-1)$.  Any quasi-exactly-solvable
operator $\tilde T_k$, can be represented by a $k$th
degree polynomial of the operators (21). If $n \leq (k-1)$, the part
of the quasi-exactly-solvable operator $\tilde T_k$ containing
derivatives up to order $n$ can be represented by a $n$th
degree polynomial in the generators (21).

(ii) Conversely, any polynomial in (21) is quasi-exactly solvable.

(iii) Among quasi-exactly-solvable operators
there exist exactly-solvable operators $\tilde E_k \subset \tilde
T_k$.}
\end{LEMMA}


\noindent
{\it Comment 3.} If we define an analogue of the universal enveloping
algebra $U_g$ for the quantum algebra $\tilde g$ as an algebra of all ordered
polynomials in generators,
then a quasi-exactly-solvable operator
$\tilde T_k$ at $k < n+1$ is simply an element of the \lq universal enveloping
algebra' $U_{sl_2({\bold R})_q}$ of the algebra $sl_2({\bold R})_q$ taken
in representation (21). If $k \geq n+1$, then $\tilde T_k$ is
represented as an element of $U_{sl_2({\bold R})_q}$ plus $B D^{n+1} $,
where $B$ is any linear difference operator of order not higher than
$(k-n-1)$.

Similar to $sl_2({\bold R})$ , one can introduce the grading of generators
(21) of $sl_2({\bold R})_q$ (see (3)) and, hence, of monomials of the
universal enveloping $U_{sl_2({\bold R})_q}$ (see (4)).
\begin{LEMMA} {\it A quasi-exactly-solvable operator
$\tilde T_k \subset U_{sl_2({\bold R})_q}$
has no terms of positive grading, iff it is an
exactly-solvable operator.}
\end{LEMMA}
\begin{THEOREM} {\it Let $n$ be a non-negative integer. Take the eigenvalue
problem for a linear difference operator of the $k$-th order in
one variable
\begin{equation}
\label{e23}
 \tilde T_k \varphi (x) \ = \ \varepsilon \varphi (x) ,
\end{equation}
where $\tilde T_k$ is symmetric. The problem (23) has $(n+1)$
linearly independent eigenfunctions in the form of a polynomial in
variable $x$ of order not higher than $n$, if and only if $T_k$ is
quasi-exactly-solvable. The problem (23) has an infinite sequence
of polynomial eigenfunctions, if and only if the operator is exactly-solvable
$\tilde E_k$.}
\end{THEOREM}

{\it Comment 4.} Saying the operator $\tilde T_k$ is symmetric, we imply
that, considering the action of this operator on a space of polynomials
of degree not higher than $n$, one can introduce a positively-defined
scalar product, and the operator $\tilde T_k$ is symmetric with respect
to it.

This theorem gives a general classification of finite-difference equations
\begin{equation}
\label{e24}
 \sum_{j=0}^{k} \tilde a_j (x) D^j \varphi (x) \ = \ \varepsilon \varphi(x)
\end{equation}
having polynomial solutions in $x$.  The coefficient functions must
have the form
\begin{equation}
\label{e25}
\tilde a_j (x) \ = \ \sum_{i=0}^{k+j} \tilde a_{j,i} x^i .
\end{equation}
In particular, this form occurs after substitution (21) into a general
$k$th degree polynomial element of the universal
enveloping algebra $U_{sl_2({\bold R})_q}$. It guarantees the existence of
at least a finite number of polynomial solutions. The coefficients
$\tilde a_{j,i}$ are related to the coefficients of the
$k$th degree polynomial element of the universal
enveloping algebra $U_{sl_2({\bold R})_q}$. The number of free parameters
of the
polynomial solutions is defined by the number of free parameters of a general
$k$-th order polynomial element of the universal
enveloping algebra $U_{sl_2({\bold R})_q}$.\footnote{For quantum
$sl_2({\bold R})_q$ algebra
there are no polynomial Casimir operators (see, e.g., \cite{z}). However,
in the representation (21) the relationship between generators analogous
to the quadratic Casimir operator
\[ q\tilde J^+\tilde J^- - \tilde J^0 \tilde J^0 + (\{ n+1 \}
- 2 \hat{n}) \tilde J^0 = \hat{n} (\hat{n} - \{ n+1 \}) \]
appears. It reduces the number of independent parameters of the
second-order polynomial element of  $U_{sl_2({\bold R})_q}$. It becomes the
standard Casimir operator at $q \rightarrow 1$. } A rather
straightforward calculation leads to the following formula
\[ par (\tilde T_k) = (k+1)^2+1 \]
(for the second-order finite-difference equation $par(\tilde T^2) = 10$).
For the case of an infinite sequence of polynomial solutions the formula
(25) simplifies to
\begin{equation}
\label{e26}
\tilde a_j (x) \ = \ \sum_{i=0}^{j} \tilde a_{j,i} x^i
\end{equation}
and the number of free parameters is given by
\[ par (\tilde E_k) = {(k+1)(k+2) \over 2} + 1 \]
(for $k=2$, $par(\tilde E^2) = 7$).
The increase in the number of free parameters
compared to ordinary differential equations is due to the presence of the
deformation parameter $q$. In \cite{t1} one can find a description
in the present approach of the $q$-deformed Hermite, Laguerre, Legendre
and Jacobi polynomials (for definitions of these polynomials see \cite{e}).
\begin{LEMMA} {\it If the operator $\tilde T_2$ (see (12)) is such that
\begin{equation}
\label{e27}
\tilde c_{++}=0 \quad and \quad \tilde c_{+} =( {\hat n}  - \{ m \})
\tilde c_{+0} \ , \ at \ some\ m=0,1,2,\dots
\end{equation}
 then the operator $\tilde T_2$ preserves both ${\cal P}_n$ and
${\cal P}_m$, and polynomial solutions in $x$ with 8 free parameters
occur.}
\end{LEMMA}

As usual in quantum algebras, a rather outstanding situation occurs
if the \linebreak deformation
parameter $q$ is equal to the root of unity. For instance, the
following statement holds.
\begin{LEMMA} {\it If a quasi-exactly-solvable operator $\tilde T_k$
preserves the space ${\cal P}_n$ and the parameter $q$ satisfies
to the equation
\begin{equation}
\label{e28}
q^n\ = \ 1 \ ,
\end{equation}
then the operator  $\tilde T_k$ preserves an infinite flag
 of polynomial spaces $ {\cal P}_0 \subset  {\cal P}_n \subset {\cal P}_{2n}
\subset \dots \subset {\cal P}_{kn} \subset \dots $.}
\end{LEMMA}

It is worth emphasizing that, in the limit as $q$ tends to one, Lemmas 2.1,2,3
and Theorem 2.1 coincide with Lemmas 1.1,2,3 and Theorem 1.1, respectively.
Thus the case of differential equations in one variable can be treated
as particular case of finite-difference ones. Evidently, one can consider
the finite-difference operators, which are a mixture of generators (21) with
the same value of $n$ and different $q$'s.
\section{Operators in one real and one Grassmann variable}
Define the following space of polynomials in $x,\theta$
\begin{equation}
\label{e29}
{\cal P}_{N,M} \ = \ \langle x^0,x^1,\dots,x^N, x^0\theta, x^1\theta, \dots,
x^M \theta \rangle
\end{equation}
where $N, M$ are non-negative integers, $x \in {\bold R}$ and $\theta$ is a
Grassmann (anticommuting) variable.

The projectivized representation of the algebra $osp(2,2)$ is given as follows.

The algebra  $osp(2,2)$ is characterized by four bosonic generators
$T^{\pm ,0},J$ and four fermionic generators $Q_{1,2},\overline{Q}_{1,2}$
and given by the commutation relations
\label{e30}
\[
[T^0 , T^{\pm}]= \pm T^{\pm} \quad , \quad [T^+ , T^-]= -2 T^0 \quad , \quad
 [J , T^{\alpha}]=0  \quad , \alpha = +,-,0 \]
\[   \{ Q_1,\overline{Q}_2 \}  = - T^- \quad , \quad \{ Q_2, \overline{Q}_1 \}
= T^+ \quad, \]
\[ {1 \over 2} (\{ \overline{Q}_1, Q_1\} + \{ \overline{Q}_2, Q_2 \})
= + J \quad ,
\quad
 {1 \over 2} (\{ \overline{Q}_1, Q_1\} - \{ \overline{Q}_2, Q_2 \}) =  T^0
\quad , \]
\[ \{ Q_1,Q_1\}=\{ Q_2,Q_2\}=\{ Q_1,Q_2\}=0\quad , \]
\[ \{ \overline{Q}_1,\overline{Q}_1\}=
\{ \overline{Q}_2,\overline{Q}_2\}=\{ \overline{Q}_1,\overline{Q}_2\}=0 \quad
,\]
\[  [Q_1 , T^+]=Q_2 \quad , \quad [Q_2 , T^+]=0 \quad , \quad [Q_1 ,
T^-]=0 \quad , \quad [Q_2 , T^-]=-Q_1 \ ,\]
\[  [\overline{Q}_1 , T^+]=0 \quad , \quad  [\overline{Q}_2 , T^+] = -
\overline{Q}_1 \quad , \quad
[\overline{Q}_1 ,T^-] = \overline{Q}_2 \quad , \quad [\overline{Q}_2 ,
T^-]=0 \ ,\]
\[    [Q_{1,2} , T^0]=\pm {1 \over 2} Q_{1,2} \quad , \quad
 [\overline{Q}_{1,2} , T^0]=\mp {1 \over 2} \overline{Q}_{1,2} \]
\begin{equation}
  [Q_{1,2}, J] = - {1 \over 2} Q_{1,2} \quad , \quad
 [ \overline{Q}_{1,2}, J] = {1 \over 2} \overline{Q}_{1,2}
\end{equation}
This algebra has the algebra $sl_2(\bold R) \oplus {\bold R}$ as sub-algebra.

The algebra (30) possesses the projectivized representation \cite{st}
\label{e31}
\[ T^+ = x^2 \partial_x - n x + x \theta \partial_{\theta},  \]
\begin{equation}
 T^0 = x \partial_x - {n \over 2} +{1 \over 2} \theta \partial _{\theta }\  ,
\end{equation}
\[ T^- = \partial_x \ .  \]
\[ J = -{ n \over 2} - {1 \over 2} \theta \partial_{\theta} \]
for bosonic (even) generators and
\begin{equation}
\label{e32}
Q= { Q_1 \brack Q_2 } = { \partial_{\theta} \brack x\partial_{\theta}} \ ,
\ {\bar Q} = { {\bar Q}_1 \brack {\bar Q}_2 } =
{ x\theta\partial_x-n \theta \brack -\theta\partial_x} ,
\end{equation}
for fermionic (odd) generators, where $x$ is a real variable and $\theta$ is
a Grassmann one. Inspection of the generators shows that if $n$ is a
non-negative integer, the representation (31), (32) is a finite-dimensional
representation of dimension $(2n+1)$. The polynomial space ${\cal P}_{n,n-1}$
describes the corresponding invariant sub-space.
\begin{definition} Let us name a linear differential operator of the $k$-th
order, ${\bf T}_k (x,\theta)$, {\it quasi-exactly-solvable} if it preserves
 the space ${\cal P}_{n,n-1}$. Correspondingly, the \linebreak operator
${\bf E}_k(x,\theta) \in {\bf T}_k(x,\theta)$, which preserves
the infinite flag $ {\cal P}_{0,0} \subset  {\cal P}_{1,0}
\subset {\cal P}_{2,1}
\subset \dots \subset {\cal P}_{n,n-1} \subset \dots $ of spaces of all
polynomials, is named {\it exactly-solvable}.
\end{definition}
\begin{LEMMA} {\it Consider the space ${\cal P}_{n,n-1}$.

 (i) Suppose $n > (k-1)$.  Any quasi-exactly-solvable operator
${\bf T}_k (x,\theta)$, can be
represented by a $k$th degree polynomial of the operators (31), (32).
 If $n \leq (k-1)$, the part of the quasi-exactly-solvable operator
${\bf T}_k (x,\theta)$ containing
derivatives in $x$ up to order $n$ can be represented by an $n$th
degree polynomial in the generators (31), (32).

(ii) Conversely, any polynomial in (31), (32) is a quasi-exactly solvable
operator.

(iii) Among quasi-exactly-solvable operators
there exist exactly-solvable operators ${\bf E}_k \subset {\bf T}_k (x,
\theta)$.}
\end{LEMMA}

Let us introduce the grading of the bosonic generators (31)
\begin{equation}
\label{e33}
deg (T^+) = +1 \ , \ deg (J,T^0) = 0 \ , \ deg (J^-) = -1
\end{equation}
and fermionic generators (32)
\begin{equation}
\label{e34}
deg (Q_2,\overline{Q}_1) =+ {1 \over 2}\ , \
deg (Q_1,\overline{Q}_2) =- {1 \over 2}
\end{equation}
Hence the grading of monomials of the generators (31), (32) is equal to
\label{e35}
\[ deg [(T^+)^{n_+} (T^0)^{n_0}(J)^{\overline{n}}(T^-)^{n_-}{Q_1}^{m_1}
{Q_2}^{m_2}
{\overline{Q}_1}^{\overline{m}_1}{\overline{Q}_2}^ {\overline{m}_2} ]
 \  = \]
\begin{equation}
  (n_+ - n_-) \ - \ (m_1 - m_2 - {\overline{m}_1} + {\overline{m}_2}) / 2
\end{equation}
The $n$'s can be arbitrary  non-negative integers, while the $m$'s are
equal to either 0 or 1. The notion of grading allows us to classify the
operators ${\bf T}_k (x,\theta)$ in the Lie-algebraic sense.
\begin{LEMMA} {\it A quasi-exactly-solvable operator
${\bf T}_k (x,\theta) \subset U_{osp(2,2)}$
 has no terms of positive grading other than  monomials of
grading +1/2 containing the generator $Q_1$ or $Q_2$,
 iff it is an exactly-solvable operator.}
\end{LEMMA}
\begin{THEOREM} {\it Let $n$ be a non-negative integer. Take the eigenvalue
problem for a linear differential
operator in one real and one Grassmann variable
\begin{equation}
\label{e36}
 {\bf T}_k (x,\theta) \varphi \ = \ \varepsilon \varphi
\end{equation}
\noindent
where $ {\bf T}_k$ is symmetric. In general,
the problem (36) has $(2n+1)$ linearly independent eigenfunctions
in the form of polynomials in variables $x,\theta$ of order
not higher than $n$, if and only if $T_k$ is quasi-exactly-solvable.
The problem (36) has an infinite sequence of polynomial eigenfunctions,
if and only if the operator is exactly-solvable.}
\end{THEOREM}

This theorem gives a general classification of differential equations
\begin{equation}
\label{e37}
 \sum_{i,j=0}^{i=k,j=1} a_{i,j} (x,\theta) \varphi^{(i,j)}_{x,\theta}
(x,\theta) \ = \ \varepsilon \varphi(x,\theta) \ ,
\end{equation}
 where the notation $ \varphi^{(i,j)}_{x,\theta}$ means the
$i$th order derivative with respect to $x$ and $j$th order derivative with
respect to $\theta$,
having at least one polynomial solution in $x,\theta$, thus resolving the
generalized Bochner problem.
Suppose that $k>0$, then the coefficient functions $a_{i,j} (x,\theta)$
should have the form
\[ a_{i,0} (x,\theta) \ =  \ \sum_{p=0}^{k+i} a_{i,0,p}\ x^p\ +\  \theta
\sum_{p=0}^{k+i-1}\overline{a}_{i,0,p}\ x^p \]
\begin{equation}
\label{e38}
a_{i,1} (x,\theta) \ = \ \sum_{p=0}^{k+i-1} a_{i,1,p}\ x^p\ + \ \theta
\sum_{p=0}^{k+i-1}
 \overline{a}_{i,1,p}\ x^p
\end{equation}
The explicit expressions (38) are obtained by substituting (31), (32) into
a general, the $k$th order, polynomial element of the universal
enveloping algebra $U_{osp(2,2)}$. Thus the coefficients $a_{i,j,p}$
can be expressed through the coefficients of the $k$-th order polynomial
element of universal enveloping algebra $U_{osp(2,2)}$. The number of free
parameters of the polynomial solutions is defined by the number of
parameters characterizing a general, $k$th order polynomial element of the
universal enveloping algebra $U_{osp(2,2)}$. However, in counting parameters
some relations between generators should be taken into account,
specific for the given representation (31), (32), like

\[ 2 T^+ J\ -\ \overline{Q}_1 Q_2\ =\ n T^+ \ , \]
\[ T^+ Q_1\ -\ T^0 Q_2\ =\ - Q_2 \ , \]
\[ T^+ \overline{Q}_2\ +\ T^0 \overline{Q} _1\ =\ {(1-n) \over 2}\
\overline{Q}_1 \ , \]
\[ J Q_2\ =\ {n \over 2}\ Q_2 \ ,\]
\[ J \overline{Q}_1\ =\ {(n+1) \over 2}\ \overline{Q}_1 \ ,\]
\[ T^+T^-\ -\ T^0 T^0\ -\ J J\ +\ T^0\ =\ - {n \over 2}(n-1) \ , \]
\[ J J\ =\ (n + {1 \over 2}) J\ -\ { n \over 4 }(n+1) \ , \]
\[ Q_1 \overline{Q}_1\ +\  Q_2 \overline{Q}_2\ -\ 2n J\ =\ - n(n+1) \ , \]
\[2 T^0 J\ +\ Q_1 \overline{Q}_1\ -\ (n+1) T^0\ -\ nJ\ =
\ - {n \over 2} (n+1) \ , \]
\[ T^- Q_2\ -\ T^0 Q_1\ =\  ( {n \over 2} +1) Q_1 \ , \]
\[ T^- \overline{Q}_1\ -\ T^0 \overline{Q}_2\ =\  {(n - 1) \over 2})
\overline{Q}_2 \ , \]
\[ J Q_1\ =\ {n \over 2} Q_1 \ ,\]
\[ J \overline{Q}_2\ =\ {n+1 \over 2} \overline{Q}_2 \ ,\]
\begin{equation}
\label{e39}
 2 J T^-\ -\ Q_1 \overline{Q}_2\ =\ (n+1) T^-
\end{equation}
between quadratic expressions in generators (and the
ideals generated by them). Straightforward analysis leads to the following
formula for the number of free \linebreak parameters
\begin{equation}
\label{e40}
par({\bf T}_k (x,\theta))= 4k(k+1)+1 \ ,
\end{equation}
 For the case of an exactly-solvable operator
(an infinite sequence of polynomial solutions of Eq. (37)), the expressions
 (38) simplify and reduce to
\[ a_{i,0} (x,\theta) \ =  \ \sum_{p=0}^{i} a_{i,0,p} x^p +  \theta
\sum_{p=0}^{i-1}\overline{a}_{i,0,p} x^p \]
\begin{equation}
\label{e41}
a_{i,1} (x,\theta) \ = \ \sum_{p=0}^{i} a_{i,1,p} x^p +
 \theta \sum_{p=0}^{i-1} \overline{a}_{i,1,p} x^p
\end{equation}
Correspondingly, the number of free parameters reduces to
\begin{equation}
\label{e42}
par({\bf E}_k (x,\theta))= 2k(k+2)+1
\end{equation}
Hence, Eq. (37) with the coefficient functions (41) gives a general form
of eigenvalue problem for the operator ${\bf T}_k$, which can lead to an
infinite family of orthogonal polynomials as eigenfunctions. If  in (41)
we put formally all coefficients, $\overline{a}_{i,0,p}$ and $a_{i,1}
(x,\theta)$ equal to zero, we reproduce the eigenvalue problem for the
differential operators in one real variable, which gives rise to all known
families of orthogonal polynomials in one real variable (see \cite{t2}).

\subsection{Second-order differential equations in $x$,$\theta$}

Now let us consider in more detail the second-order differential equation
Eq. (37), which can possess polynomial solutions. As follows from Theorem 3.1,
the corresponding differential operator ${\bf T}_2(x,\theta)$ should be
quasi-exactly-solvable. Hence, this operator can
be expressed in terms of $osp(2,2)$ generators taking into account the
relations (39)
\label{e43}
\[ {\bf T}_2 =  c_{++} T^+ T^+ + c_{+0} T^+ T^0 +  c_{+-} T^+ T^- +
 c_{0-} T^0 T^- + c_{--} T^- T^- + \]
\[ c_{+J} T^+ J + c_{0J} T^0 J + c_{-J} T^- J + \]
\[ c_{+ \overline{1}} T^+ \overline{Q}_1 +  c_{+2} T^+ Q_2 +
  c_{+1} T^+ Q_1 + c_{+ \overline{2}} T^+ \overline{Q}_2 +
   c_{01} T^0 Q_1 + \]
\[ c_{0 \overline{2}} T^0 \overline{Q}_2 +  c_{-1} T^- Q_1 +
c_{- \overline{2}} T^- \overline{Q}_2 + \]
\begin{equation}
c_+ T^+ + c_0 T^0 + c_- T^- + c_J J + c_1 Q^1 + c_2 Q^2 +  c_{\overline{1}}
\overline{Q}_1 + c_{\overline{2}} \overline{Q}_2 + c
\end{equation}
where $c_{\alpha \beta}, c_{\alpha}, c $ are parameters.
Following Lemma 3.2, under the conditions
\begin{equation}
\label{e44}
 c_{++}  = c_{+0} = c_{+ \overline{1}}= c_{+ \overline{2}}=
c_{ \overline{1}} = c_{+2}  = c_{+J}= c_{+}  =0 \ ,
\end{equation}
the operator ${\bf T}_2(x,\theta)$ becomes exactly-solvable (see (41)).

Now we proceed to the detailed analysis of the quasi-exactly-solvable
operator ${\bf T}_2(x,\theta)$. Set
\begin{equation}
\label{e45}
 c_{++}  = 0
\end{equation}
in Eq. (43). The remainder possesses an exceptionally rich structure. The
whole situation can be subdivided into three cases
\begin{equation}
\label{e46}
 c_{+2} \neq 0 \ , \  c_{+ \overline{1}}\ =\ 0 \ (case\ I)
\end{equation}
\begin{equation}
\label{e47}
 c_{+2}\ =\ 0 \ , \  c_{+ \overline{1}} \neq 0 \ (case \ II)
\end{equation}
\begin{equation}
\label{e48}
 c_{+2}\ =\ 0 \ , \  c_{+ \overline{1}}\ =\ 0 \ (case \ III)
\end{equation}
We emphasize that we keep the parameter $n$ in the representation (31), (32)
as a fixed, non-negative integer.
\newpage
{\bf Case I.} The conditions (45) and (46) are fulfilled (see Fig. 3.I).
\vskip 0.6truecm
{\it Case I.1.1.}  If
\[ (n+2) c_{+0} + n c_{+J} + 2 c_{+} = 0\ ,\]
\[ c_{+\overline{2}}=c_{\overline{1}}=0\ , \]
\begin{equation}
\label{e49}
(n+1) c_{0\overline{2}}+2 c_{\overline{2}}=0\ ,
\end{equation}
 then ${\bf T}_2(x,\theta)$ preserves ${\cal P}_{n,n-1}$ and
{\bf ${\cal P}_{n+1,n-1}$}.

{\it Case I.1.2.}  If
\[ (n+4+2m) c_{+0} + n c_{+J} + 2 c_{+} = 0\ ,\]
\[ c_{+\overline{2}}= c_{\overline{1}}=0\ , \]
\begin{equation}
\label{e50}
 c_{0\overline{2}}= c_{\overline{2}}= c_{-\overline{2}} = 0\ ,
\end{equation}
at a certain integer $m \geq 0$ , then ${\bf T}_2(x,\theta)$ preserves
${\cal P}_{n,n-1}$ and {\bf ${\cal P}_{n+2+m,n-1}$}. If $m$ is non-integer,
then ${\bf T}_2(x,\theta)$
preserves ${\cal P}_{n,n-1}$ and {\bf ${\cal P}_{\infty,n-1}$}.

{\it Case I.1.3.}  If
\[ (n+1) c_{+J} + 2 c_{+} = 0\ ,\]
\[ c_{+0}=0\ , \]
\[ c_{+\overline{2}}= c_{\overline{1}}=0\ , \]
\begin{equation}
\label{e51}
 c_{0\overline{2}}= c_{\overline{2}}= c_{-\overline{2}} = 0\ ,
\end{equation}
 then ${\bf T}_2(x,\theta)$ preserves the infinite flag
of polynomial spaces the {\bf ${\cal P}_{n+m,n-1}, \linebreak
m=0,1,2,\ldots$}.

{\it Case I.2.1.}  If
\[ (n-3) c_{+0} + (n+1) c_{+J} + 2 c_{+} = 0\ ,\]
\[ (n-1)c_{+\overline{2}}=c_{\overline{1}}\ , \]
\begin{equation}
\label{e52}
(n-1) c_{0\overline{2}}+2 c_{\overline{2}}=0\ ,
\end{equation}
 then ${\bf T}_2(x,\theta)$ preserves ${\cal P}_{n,n-1}$ and
{\bf ${\cal P}_{n,n-2}$}.

{\it Case I.2.2.}  If
\[ 3 c_{+0} - c_{+J}  = 0\ ,\]
\[ (2k+2n+4) c_{+0} + 2 c_{+} = 0\ ,\]
\[ c_{+\overline{2}}=c_{\overline{1}}=0\ , \]
\begin{equation}
\label{e53}
(2k-n+3) c_{0\overline{2}}+2 c_{\overline{2}}=0\ ,
\end{equation}
at a certain integer $k \geq 0$, then ${\bf T}_2(x,\theta)$ preserves
${\cal P}_{n,n-1}$ and {\bf ${\cal P}_{k+2,k}$}.

{\it Case I.2.3.}  If
\[ c_{+0} = c_{+J}= c_{+} = 0\ ,\]
\[ c_{+\overline{2}}=c_{\overline{1}}=0\ , \]
\begin{equation}
\label{e54}
 c_{0\overline{2}}= c_{\overline{2}}=0\ ,
\end{equation}
 then ${\bf T}_2(x,\theta)$ preserves ${\cal P}_{n,n-1}$ and the infinite
flag of the polynomial spaces {\bf ${\cal P}_{k+2,k}, \linebreak
 k=0,1,2,\ldots$}.
Note in general for this case $ c_{-\overline{2}} \neq 0$.

{\it Case I.3.1.}  If
\[ (n-5-2m) c_{+0} + (n+1) c_{+J} + 2 c_{+} = 0\ ,\]
\[ c_{+\overline{2}}= c_{\overline{1}}=0\ , \]
\begin{equation}
\label{e55}
 c_{0\overline{2}}= c_{\overline{2}}= c_{-\overline{2}} = 0\ ,
\end{equation}
at a certain integer $0 \leq m \leq (n-3)$ , then ${\bf T}_2(x,\theta)$
preserves ${\cal P}_{n,n-1}$ and {\bf ${\cal P}_{n,n-3-m}$}.

{\it Case I.3.2.}  If
\[ c_{+0}=0\ , \]
\[ (n+1) c_{+J} + 2 c_{+} = 0\ ,\]
\[ c_{+\overline{2}}= c_{\overline{1}}=0\ , \]
\begin{equation}
\label{e56}
 c_{0\overline{2}}= c_{\overline{2}}= c_{-\overline{2}} = 0\ ,
\end{equation}
 then ${\bf T}_2(x,\theta)$ preserves ${\cal P}_{n,n-1}$ and the sequence
of the polynomial spaces {\bf ${\cal P}_{n,n-3-m},\linebreak m=0,1,2,\ldots ,
 (n-3)$}.

{\it Case I.3.3.}  If
\[ (2k+1-n) c_{+0} + (n+1) c_{+J} + 2 c_{+} = 0\ ,\]
\[ (2m+5) c_{+0} - c_{+J} = 0\ ,\]
\[ c_{+\overline{2}}= c_{\overline{1}}=0\ , \]
\begin{equation}
\label{e57}
 c_{0\overline{2}}= c_{\overline{2}}= c_{-\overline{2}} = 0\ ,
\end{equation}
at  certain integers $k , m \geq 0 $ , then ${\bf T}_2(x,\theta)$ preserves
${\cal P}_{n,n-1}$ and {\bf ${\cal P}_{k+3+m,k}$}.

{\it Case I.3.4.}  If
\[  c_{+0} = c_{+J}= c_{+} = 0\ ,\]
\[ c_{+\overline{2}}= c_{\overline{1}}=0\ , \]
\begin{equation}
\label{e58}
 c_{0\overline{2}}= c_{\overline{2}}= c_{-\overline{2}} = 0\ ,
\end{equation}
 then ${\bf T}_2(x,\theta)$ preserves
${\cal P}_{n,n-1}$ and the infinite flag of polynomial spaces {\bf ${\cal
P}_{k+3+m,k},\linebreak  k,m=0,1,2,\ldots$} (cf. {\it Cases I.1.3 and I.2.3}).
\vskip 0.6truecm
{\bf Case II.} The conditions (45) and (47) are fulfilled (see Fig. 3.II).
\vskip 0.6truecm
{\it Case II.1.1.}  If
\[ (n+1) c_{+0} + (n+1) c_{+J} + 2 c_{+} = 0\ ,\]
\begin{equation}
\label{e59}
 c_{2}= 0 ,
\end{equation}
 then ${\bf T}_2(x,\theta)$ preserves ${\cal P}_{n,n-1}$ and
{\bf ${\cal P}_{n,n}$}.

{\it Case II.1.2.}  If
\[ (n+3) c_{+0} + (n+1) c_{+J} + 2 c_{+} = 0\ ,\]
\[(n+2) c_{01}+ 2c_{1}=0\ , \]
\begin{equation}
\label{e60}
 c_{+1} = c_{2} = 0\ ,
\end{equation}
  then ${\bf T}_2(x,\theta)$ preserves ${\cal P}_{n,n-1}$ and
{\bf ${\cal P}_{n,n+1}$}.

{\it Case II.1.3.}  If
\[ (2k+5+n) c_{+0} + (n+1) c_{+J} + 2 c_{+} = 0\ ,\]
\[  c_{01} = c_{1} = 0\ ,\]
\[ c_{+1} = c_{2} = 0\ , \]
\begin{equation}
\label{e61}
 c_{-1} = 0\ ,
\end{equation}
at a certain integer $k \geq 0$ , then ${\bf T}_2(x,\theta)$ preserves
${\cal P}_{n,n-1}$ and {\bf ${\cal P}_{n,n+2+k}$}.

{\it Case II.1.4.}  If
\[  c_{+0} = 0 , \]
\[ (n+1) c_{+J}+2 c_{+} = 0\ ,\]
\[  c_{01} = c_{1} = 0\ ,\]
\[ c_{+1} = c_{2} = 0\ , \]
\begin{equation}
\label{e62}
 c_{-1} = 0 ,
\end{equation}
 then ${\bf T}_2(x,\theta)$ preserves
${\cal P}_{n,n-1}$ and the infinite flag of polynomial spaces {\bf ${\cal
P}_{n,n+k} , \linebreak  k = 0,1,2,\ldots$}
 (cf. {\it Cases I.1.3, I.2.3 and I.3.4}).

{\it Case II.2.1.}  If
\[ (n-2) c_{+0} + n c_{+J} + 2 c_{+} = 0\ ,\]
\begin{equation}
\label{e63}
 c_{+1} = c_{2} \ ,
\end{equation}
then ${\bf T}_2(x,\theta)$ preserves ${\cal P}_{n,n-1}$ and
{\bf ${\cal P}_{n-1,n-1}$}.

{\it Case II.2.2.}  If
\[ (n+1) c_{+0} -  c_{+} = 0\ ,\]
\[ 3 c_{+0} +  c_{+J}  = 0\ ,\]
\[ n c_{01}+ 2c_{1}=0\ , \]
\begin{equation}
\label{e64}
 c_{+1} = c_{2}=0\ ,
\end{equation}
then ${\bf T}_2(x,\theta)$ preserves ${\cal P}_{n,n-1}$ and
{\bf ${\cal P}_{n-1,n}$}.

{\it Case II.2.3.}  If
\[ (n-2) c_{+0} + n c_{+J} + 2 c_{+} = 0\ ,\]
\[ (2k+1) c_{+0} +  c_{+J} = 0\ ,\]
\[ c_{+1} = c_{2} = 0\ , \]
\begin{equation}
\label{e65}
 c_{01} = c_{1} = c_{-1}= 0\ ,
\end{equation}
at a  certain integer $k \geq 0$,  then ${\bf T}_2(x,\theta)$
preserves ${\cal P}_{n,n-1}$ and {\bf ${\cal P}_{n-1,n+k+1}$}.

{\it Case II.2.4.}  If
\[  c_{+0} = c_{+J}= c_{+} = 0\ ,\]
\[ c_{+1} = c_{2} = 0\ , \]
\begin{equation}
\label{e66}
 c_{01} = c_{1} = c_{-1}= 0\ ,
\end{equation}
 then ${\bf T}_2(x,\theta)$ preserves
${\cal P}_{n,n-1}$ and the infinite flag of the polynomial spaces
{\bf ${\cal P}_{n-1,n+k},\linebreak  k = 0,1,2,\ldots$}
 (cf. {\it Cases I.1.3, I.2.3, I.3.4 and II.1.4}).

{\it Case II.3.1.}  If
\[ (n-4)c_{+0}  + nc_{+J} + 2c_{+} = 0\ ,\]
\[ (n-2) c_{01}+ 2c_{1}=0\ , \]
\begin{equation}
\label{e67}
 c_{+1} = c_{2} = 0\ ,
\end{equation}
 then ${\bf T}_2(x,\theta)$ preserves ${\cal P}_{n,n-1}$ and
{\bf ${\cal P}_{n-2,n-1}$}.

{\it Case II.4.1.}  If
\[ (m-2n) c_{+0} +  c_{+} = 0\ ,\]
\[ 3 c_{+0} +  c_{+J}  = 0\ ,\]
\[ (2m+2-n) c_{01}+ 2c_{1}=0\ , \]
\begin{equation}
\label{e68}
 c_{+1} = c_{2} = 0\ ,
\end{equation}
at a  certain integer $m \geq 0$, then ${\bf T}_2(x,\theta)$ preserves
${\cal P}_{n,n-1}$ and {\bf ${\cal P}_{m,m+1}$}.

{\it Case II.4.2.}  If
\[ (2m-n)c_{+0}  + nc_{+J} + 2c_{+} = 0\ ,\]
\[ (2k+5) c_{+0} + 2 c_{+J}  = 0\ ,\]
\[ c_{+1} = c_{2}=0 , \]
\begin{equation}
\label{e69}
 c_{01} = c_{1} = c_{-1}= 0
\end{equation}
at  certain integers $k \geq 0 , m \geq 0$, then ${\bf T}_2(x,\theta)$
preserves ${\cal P}_{n,n-1}$ and {\bf ${\cal P}_{m,m+2+k}$}.

{\it Case II.4.3.}  If
\[ c_{+0}  = c_{+J} = c_{+} = 0\ ,\]
\[ c_{+1} = c_{2}=0 , \]
\begin{equation}
\label{e70}
 c_{01} = c_{1} = c_{-1}= 0
\end{equation}
 then ${\bf T}_2(x,\theta)$ preserves ${\cal P}_{n,n-1}$ and the infinite
flag of the polynomial spaces \linebreak
{\bf ${\cal P}_{m,m+1+k},\ m,k=0,1,2,\ldots $}
 (cf. {\it Cases I.1.3, I.2.3, I.3.4, II.1.4 and II.2.4}).
\vskip 0.6truecm
{\bf Case III.} The conditions (45) and (48) are fulfilled (see Fig. 3.III).
\vskip 0.6truecm
{\it Case III.1.1.}  If
\[ (2m-n)c_{+0}  + nc_{+J} + 2c_{+} = 0\ ,\]
\[  c_{+0} +  c_{+J}  = 0\ ,\]
\begin{equation}
\label{e71}
 (m-n) c_{+1} + c_{2}=0 ,
\end{equation}
at a certain integer $ m \geq 0$, then ${\bf T}_2(x,\theta)$ preserves
${\cal P}_{n,n-1}$ and {\bf ${\cal P}_{m,m}$}.

\newpage
{\it Case III.1.2.}  If
\[ c_{+0}  = c_{+J} = c_{+} = 0\ ,\]
\begin{equation}
\label{e72}
 c_{+1} = c_{2}=0 ,
\end{equation}
 then ${\bf T}_2(x,\theta)$ preserves ${\cal P}_{n,n-1}$ and the infinite
flag of the polynomial spaces {\bf ${\cal P}_{m,m},\linebreak  m=0,1,2,\ldots
 $}
 (cf. {\it Cases I.1.3, I.2.3, I.3.4, II.1.4, II.2.4 and
 II.4.3}).

{\it Case III.2.1.}  If
\[ (2m-n)c_{+0}  + nc_{+J} + 2c_{+} = 0\ ,\]
\[  c_{+0} - c_{+J}  = 0\ ,\]
\begin{equation}
\label{e73}
 m c_{+ \overline{2}} - c_{\overline{1}}=0 ,
\end{equation}
at a certain integer $ m \geq 0$, then ${\bf T}_2(x,\theta)$ preserves
${\cal P}_{n,n-1}$ and {\bf ${\cal P}_{m,m-1}$}.

{\it Case III.2.2.}  If
\[ c_{+0}  = c_{+J} = c_{+} = 0\ ,\]
\begin{equation}
\label{e74}
 c_{+ \overline{2}} = c_{\overline{1}}=0
\end{equation}
 then ${\bold T}_2(x,\theta)$ preserves ${\cal P}_{n,n-1}$ and the infinite
flag of polynomial spaces {\bf ${\cal P}_{m,m-1},\linebreak  m=0,1,2,\ldots $}
 (cf. {\it Cases I.1.3, I.2.3, I.3.4, II.1.4, II.2.4,
II.4.3 and III.1.2}). This case corresponds to exactly-solvable operators
${\bf E}_k$.
\vskip 0.6truecm

In \cite{t2} it has been shown that under a certain condition some
quasi-exactly-solvable operators $T_2(x)$ in one real
variable can preserve two polynomial spaces of different dimensions
$n$ and $m$ (see Lemma 1.3). It has been shown that
those quasi-exactly-solvable operators  $T_2(x)$ can be represented through
the generators of $sl_2({\bold R})$ in a projectivized representation
characterized either by the mark $n$ or by the mark $m$.
The above analysis shows that the quasi-exactly-solvable operators
${\bf T}_2(x,\theta)$
in two variables (one real and one Grassmann) possess an extremely rich
variety of internal properties. They are characterized by different
structures of invariant sub-spaces. However, generically the
quasi-exactly-solvable operators  ${\bf T}_2(x,\theta)$ can preserve
either one, or two,
or infinitely many polynomial spaces. For the latter, those operators become
\lq exactly-solvable' (see  {\it Cases I.1.3, I.2.3, I.3.4,
 II.1.4, II.2.4, II.4.3 and III.1.2})
  giving rise to eigenvalue problems (36) possessing
infinite sequences of polynomial eigenfunctions. In general, for the two
latter cases the interpretation of  ${\bf T}_2(x,\theta)$
in terms of $osp(2,2)$ generators characterized by \linebreak different marks
 does not exist, unlike the case of quasi-exactly-solvable operators
in one real variable. The only exceptions are given  by
{\it Case III.2.1} and {\it Case III.2.2}.
\vskip 1truecm
\newpage
\begin{center}
FIGURES. 3.0--3.III
\end{center}
\vskip 1truecm

Newton diagrams describing invariant subspaces ${\cal P}_{N,M}$
 of the second-order polynomials in the generators of $osp(2,2)$.
The lower line corresponds to the part of the space of zero
degree in  $\theta$ and the upper line of first degree in $\theta$.
The letters without brackets indicate the maximal degree of the polynomial
in $x$. The letters in brackets indicate the maximal (or minimal) possible
degree, if the degree can be varied. The {\it thin} line displays
schematically the length of polynomial in $x$ (the number of monomials).
The {\it thick} line shows that the length of polynomial can not be more
(or less) than that size. The {\it dashed} line means that the length of
polynomial can take any size on this line. If the dashed line is unbounded,
it means that the degree of the polynomial can be arbitrary up to infinity.
The numbering of the figures I-III corresponds to the cases, which satisfy the
conditions (45), (46) ({\it Case I}); (45), (47) ({\it Case II}) and (45),
(48) ({\it Case III}).

\vskip 1.5truecm
\noindent
\unitlength.8pt
\begin{picture}(400,50)(-10,-20)
\linethickness{1.2pt}
\put(10,10){\circle*{5}}
\put(10,10){\line(1,0){50}}
\put(60,10){\circle*{5}}
\put(10,20){\circle*{5}}
\put(10,20){\line(1,0){40}}
\put(50,20){\circle*{5}}
\put(-10,7){$x$}
\put(-10,17){$\Theta$}
\put(35,25){$n-1$}
\put(60,0){$n$}
\end{picture}
\begin{center}
Fig. 3.0. \ Basic subspace
\end{center}
\vskip 0.3truecm
\begin{picture}(400,50)(-10,-20)
\linethickness{1.2pt}
\put(10,10){\circle*{5}}
\put(10,10){\line(1,0){50}}
\put(60,10){\circle*{5}}
\put(10,20){\circle*{5}}
\put(10,20){\line(1,0){30}}
\put(40,20){\circle*{5}}
\put(-10,7){$x$}
\put(-10,17){$\Theta$}
\put(25,25){$n-1$}
\put(50,-2){$n+1$}
\put(22,-20){(a)}
\put(130,10){\circle*{5}}
\linethickness{2.4pt}
\put(130,10){\line(1,0){40}}
\put(170,10){\circle*{5}}
\linethickness{1.2pt}
\put(170,10){\line(1,0){30}}
\put(200,10){\circle*{5}}
\put(130,20){\circle*{5}}
\put(130,20){\line(1,0){30}}
\put(160,20){\circle*{5}}
\put(150,25){$n-1$}
\put(148,-2){$(n+2)$}
\put(192,-2){$n+2+m$}
\put(182,-20){(b)}
\linethickness{2.4pt}
\put(290,10){\circle*{5}}
\put(290,10){\line(1,0){40}}
\put(330,10){\circle*{5}}
\linethickness{1.2pt}
\put(290,20){\circle*{5}}
\put(290,20){\line(1,0){30}}
\put(320,20){\circle*{5}}
\put(320,10){\dashbox{2}(40,0){}}
\put(310,25){$n-1$}
\put(315,-2){$(n+2)$}
\put(320,-20){(c)}
\end{picture}
\begin{center}
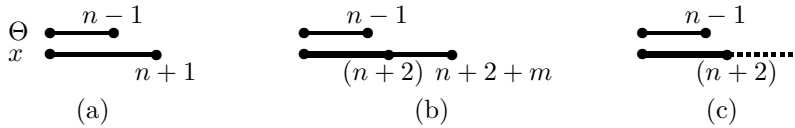
Fig. 3.I.1 . \  Subspaces for the {\it Case I.1}
\end{center}
\vskip 0.3truecm
\begin{picture}(400,50)(-10,-20)
\linethickness{1.2pt}
\put(10,10){\circle*{5}}
\put(10,10){\line(1,0){50}}
\put(60,10){\circle*{5}}
\put(10,20){\circle*{5}}
\put(10,20){\line(1,0){25}}
\put(35,20){\circle*{5}}
\put(-10,7){$x$}
\put(-10,17){$\Theta$}
\put(25,25){$n-2$}
\put(57,-2){$n$}
\put(22,-20){(a)}
\put(130,10){\circle*{5}}
\linethickness{2.4pt}
\put(130,10){\line(1,0){10}}
\put(140,10){\circle*{5}}
\linethickness{1.2pt}
\put(140,10){\line(1,0){60}}
\put(200,10){\circle*{5}}
\put(130,20){\circle*{5}}
\put(130,20){\line(1,0){30}}
\put(160,20){\circle*{5}}
\put(157,25){$k$}
\put(133,-3){$(2)$}
\put(192,-2){$k+2$}
\put(182,-20){(b)}
\linethickness{2.4pt}
\put(290,10){\circle*{5}}
\put(290,10){\line(1,0){10}}
\put(300,10){\circle*{5}}
\linethickness{1.2pt}
\put(290,20){\circle*{5}}
\put(290,20){\dashbox{2}(40,0){}}
\put(300,10){\dashbox{2}(60,0){}}
\put(295,-3){$(2)$}
\put(327,25){$k$}
\put(352,-2){$k+2$}
\put(320,-20){(c)}
\end{picture}
\begin{center}
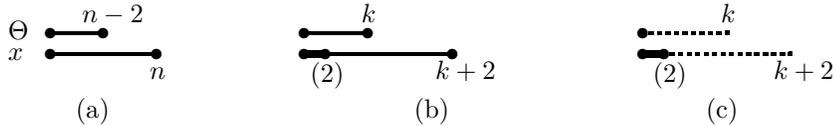
Fig. 3.I.2 . \  Subspaces for the {\it Case I.2}
\end{center}
\vskip 0.3truecm
\begin{picture}(400,50)(-10,-20)
\linethickness{1.2pt}
\put(10,10){\circle*{5}}
\put(10,10){\line(1,0){50}}
\put(60,10){\circle*{5}}
\linethickness{2.4pt}
\put(10,20){\circle*{5}}
\put(10,20){\line(1,0){40}}
\put(25,20){\circle*{7}}
\put(5,28){$n-3+m$}
\put(50,20){\circle*{5}}
\put(-10,7){$x$}
\put(-10,17){$\Theta$}
\put(60,21){$(n-3)$}
\put(60,0){$n$}
\put(22,-20){(a)}
\linethickness{1.2pt}
\put(120,10){\circle*{5}}
\put(120,10){\line(1,0){60}}
\put(180,10){\circle*{5}}
\put(120,20){\circle*{5}}
\put(120,20){\dashbox{2}(40,0){}}
\put(160,20){\circle*{5}}
\put(165,21){$(n-3)$}
\put(180,0){$n$}
\put(142,-20){(b)}
\put(250,10){\circle*{5}}
\linethickness{2.4pt}
\put(250,10){\line(1,0){12}}
\put(262,10){\circle*{5}}
\linethickness{1.2pt}
\put(262,10){\line(1,0){50}}
\put(312,10){\circle*{5}}
\put(250,20){\circle*{5}}
\put(250,20){\line(1,0){30}}
\put(280,20){\circle*{5}}
\put(277,25){$k$}
\put(256,-3){$(3)$}
\put(282,-2){$k+3+m$}
\put(280,-20){(c)}
\end{picture}
\vskip 0.3truecm
\begin{picture}(400,50)(-10,-20)
\linethickness{1.2pt}
\put(10,10){\circle*{5}}
\linethickness{2.4pt}
\put(10,10){\line(1,0){12}}
\put(22,10){\circle*{5}}
\linethickness{1.2pt}
\put(22,10){\dashbox{2}(50,0){}}
\put(10,20){\circle*{5}}
\put(10,20){\dashbox{2}(30,0){}}
\put(37,25){$k$}
\put(16,-3){$(3)$}
\put(42,-2){$k+3+m$}
\put(40,-20){(d)}
\end{picture}
\begin{center}
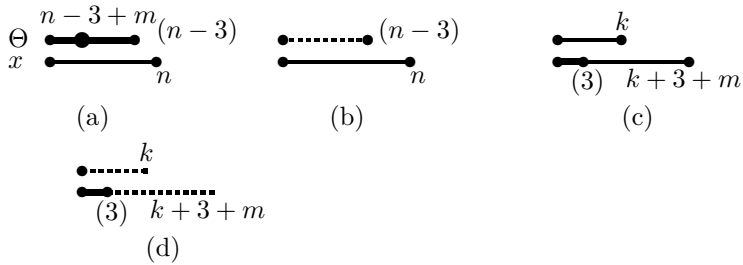
Fig. 3.I.3 . \  Subspaces for the {\it Case I.3}
\end{center}
\vskip 0.3truecm
\begin{picture}(400,50)(-10,-20)
\linethickness{1.2pt}
\put(10,10){\circle*{5}}
\put(10,10){\line(1,0){30}}
\put(40,10){\circle*{5}}
\put(10,20){\circle*{5}}
\put(10,20){\line(1,0){30}}
\put(40,20){\circle*{5}}
\put(-10,7){$x$}
\put(-10,17){$\Theta$}
\put(45,25){$n$}
\put(45,0){$n$}
\put(20,-20){(a)}
\put(110,10){\circle*{5}}
\put(110,10){\line(1,0){30}}
\put(140,10){\circle*{5}}
\put(110,20){\circle*{5}}
\put(110,20){\line(1,0){40}}
\put(150,20){\circle*{5}}
\put(145,25){$n+1$}
\put(145,0){$n$}
\put(120,-20){(b)}
\put(230,10){\circle*{5}}
\put(230,20){\circle*{5}}
\linethickness{2.4pt}
\put(230,20){\line(1,0){40}}
\put(270,20){\circle*{5}}
\linethickness{1.2pt}
\put(270,20){\line(1,0){20}}
\put(290,20){\circle*{5}}
\put(230,10){\line(1,0){25}}
\put(255,10){\circle*{5}}
\put(245,28){$(n+2)$}
\put(289,25){$n+2+k$}
\put(250,-2){$n-1$}
\put(252,-20){(c)}
\end{picture}
\vskip 0.3truecm
\begin{picture}(400,50)(-10,-20)
\linethickness{1.2pt}
\put(10,10){\circle*{5}}
\put(10,20){\circle*{5}}
\linethickness{2.4pt}
\put(10,20){\line(1,0){30}}
\put(40,20){\circle*{5}}
\linethickness{1.2pt}
\put(40,20){\dashbox{2}(35,0){}}
\put(10,10){\line(1,0){30}}
\put(40,10){\circle*{5}}
\put(35,28){$(n)$}
\put(71,25){$n+k$}
\put(38,-2){$n$}
\put(40,-20){(d)}
\end{picture}
\begin{center}
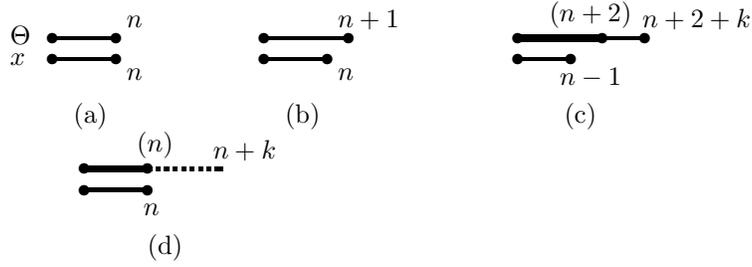
Fig. 3.II.1 . \ Subspaces for the {\it Case II.1}
\end{center}
\vskip 0.3truecm
\noindent
\begin{picture}(400,50)(-10,-20)
\linethickness{1.2pt}
\put(10,10){\circle*{5}}
\put(10,10){\line(1,0){30}}
\put(40,10){\circle*{5}}
\put(10,20){\circle*{5}}
\put(10,20){\line(1,0){30}}
\put(40,20){\circle*{5}}
\put(-10,7){$x$}
\put(-10,17){$\Theta$}
\put(43,25){$n-1$}
\put(43,0){$n-1$}
\put(20,-20){(a)}
\put(110,10){\circle*{5}}
\put(110,10){\line(1,0){30}}
\put(140,10){\circle*{5}}
\put(110,20){\circle*{5}}
\put(110,20){\line(1,0){40}}
\put(150,20){\circle*{5}}
\put(145,25){$n$}
\put(138,0){$n-1$}
\put(120,-20){(b)}
\put(230,10){\circle*{5}}
\put(230,20){\circle*{5}}
\linethickness{2.4pt}
\put(230,20){\line(1,0){35}}
\put(265,20){\circle*{5}}
\linethickness{1.2pt}
\put(265,20){\line(1,0){35}}
\put(300,20){\circle*{5}}
\put(230,10){\line(1,0){25}}
\put(255,10){\circle*{5}}
\put(245,28){$(n+1)$}
\put(289,25){$n+1+k$}
\put(250,-2){$n-1$}
\put(252,-20){(c)}
\end{picture}
\vskip 0.3truecm
\begin{picture}(400,50)(-10,-20)
\linethickness{1.2pt}
\put(10,10){\circle*{5}}
\put(10,20){\circle*{5}}
\linethickness{2.4pt}
\put(10,20){\line(1,0){35}}
\put(45,20){\circle*{5}}
\linethickness{1.2pt}
\put(45,20){\dashbox{2}(35,0){}}
\put(10,10){\line(1,0){25}}
\put(35,10){\circle*{5}}
\put(35,28){$(n)$}
\put(71,25){$n+k$}
\put(38,-2){$n-1$}
\put(40,-20){(d)}
\end{picture}
\begin{center}
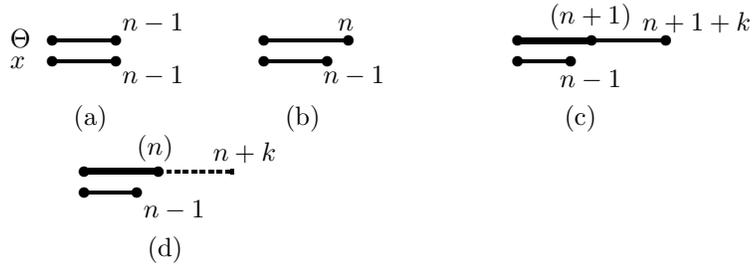
Fig. 3.II.2 . \ Subspaces for the {\it Case II.2}
\end{center}
\vskip 0.3truecm
\begin{picture}(400,50)(-10,-20)
\linethickness{1.2pt}
\put(10,10){\circle*{5}}
\put(10,10){\line(1,0){30}}
\put(40,10){\circle*{5}}
\put(10,20){\circle*{5}}
\put(10,20){\line(1,0){40}}
\put(50,20){\circle*{5}}
\put(-10,7){$x$}
\put(-10,17){$\Theta$}
\put(53,25){$n-1$}
\put(43,0){$n-2$}
\end{picture}
\begin{center}
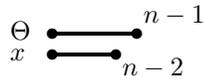
Fig. 3.II.3 . \ Subspace for the {\it Case II.3}
\end{center}
\vskip 0.3truecm
\begin{picture}(400,50)(-10,-20)
\put(10,20){\circle*{5}}
\linethickness{2.4pt}
\put(10,20){\line(1,0){10}}
\put(20,20){\circle*{5}}
\linethickness{1.2pt}
\put(20,20){\line(1,0){40}}
\put(60,20){\circle*{5}}
\put(10,10){\circle*{5}}
\put(10,10){\line(1,0){40}}
\put(50,10){\circle*{5}}
\put(57,25){$m+1$}
\put(16,28){$(1)$}
\put(52,-2){$m$}
\put(22,-20){(a)}
\put(130,20){\circle*{5}}
\linethickness{2.4pt}
\put(130,20){\line(1,0){20}}
\put(150,20){\circle*{5}}
\linethickness{1.2pt}
\put(150,20){\line(1,0){40}}
\put(190,20){\circle*{5}}
\put(130,10){\circle*{5}}
\put(130,10){\line(1,0){40}}
\put(170,10){\circle*{5}}
\put(185,25){$m+2+k$}
\put(145,28){$(2)$}
\put(170,-2){$m$}
\put(142,-20){(b)}
\put(250,20){\circle*{5}}
\linethickness{2.4pt}
\put(250,20){\line(1,0){20}}
\put(270,20){\circle*{5}}
\linethickness{1.2pt}
\put(270,20){\dashbox{2}(40,0){}}
\put(250,10){\circle*{5}}
\put(250,10){\dashbox{2}(40,0){}}
\put(305,25){$m+2+k$}
\put(265,28){$(2)$}
\put(290,-2){$m$}
\put(272,-20){(c)}
\end{picture}
\begin{center}
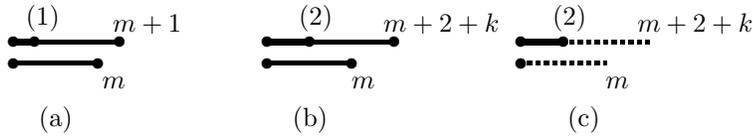
Fig. 3.II.4 . \ Subspaces for the {\it Case II.4}
\end{center}
\vskip 0.3truecm
\begin{picture}(400,50)(-10,-20)
\linethickness{1.2pt}
\put(10,10){\circle*{5}}
\put(10,10){\line(1,0){30}}
\put(40,10){\circle*{5}}
\put(10,20){\circle*{5}}
\put(10,20){\line(1,0){30}}
\put(40,20){\circle*{5}}
\put(-10,7){$x$}
\put(-10,17){$\Theta$}
\put(43,25){$m$}
\put(43,0){$m$}
\put(20,-20){(a)}
\linethickness{1.2pt}
\put(120,10){\circle*{5}}
\put(120,10){\dashbox{2}(40,0){}}
\put(120,20){\circle*{5}}
\put(120,20){\dashbox{2}(40,0){}}
\put(163,25){$m$}
\put(163,0){$m$}
\put(130,-20){(b)}
\end{picture}
\begin{center}
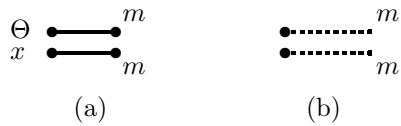
Fig. 3.III.1 . \ Subspaces for the {\it Case III.1}
\end{center}
\vskip 0.3truecm
\begin{picture}(400,50)(-10,-20)
\put(10,10){\circle*{5}}
\linethickness{2.4pt}
\put(10,10){\line(1,0){10}}
\put(20,10){\circle*{5}}
\linethickness{1.2pt}
\put(20,10){\line(1,0){50}}
\put(70,10){\circle*{5}}
\put(10,20){\circle*{5}}
\put(10,20){\line(1,0){50}}
\put(60,20){\circle*{5}}
\put(62,25){$m-1$}
\put(17,-3){$(1)$}
\put(72,-2){$m$}
\put(42,-20){(a)}
\put(120,10){\circle*{5}}
\linethickness{2.4pt}
\put(120,10){\line(1,0){10}}
\put(130,10){\circle*{5}}
\linethickness{1.2pt}
\put(130,10){\dashbox{2}(50,0){}}
\put(120,20){\circle*{5}}
\put(120,20){\dashbox{2}(50,0){}}
\put(172,25){$m-1$}
\put(127,-3){$(1)$}
\put(182,-2){$m$}
\put(152,-20){(b)}
\end{picture}
\begin{center}
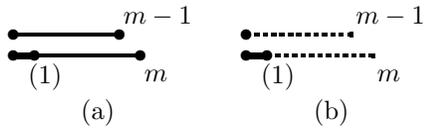
Fig. 3.III.2 . \ Subspaces for the {\it Case III.2}
\end{center}
\newpage
\subsection{2 x 2 matrix differential equations in $x$}
It is well known that anti-commuting variables can be represented by
matrices.
In our case the matrix representation is as follows: substitute $\theta$ and
$\partial_{\theta}$ in the generators (31), (32) by the Pauli matrices
$\sigma^+$ and $\sigma^-$, respectively, acting on two-component spinors.
In fact, all main notations are preserved such as quasi-exactly-solvable and
exactly-solvable operator, grading etc.

In the explicit form the fermionic generators (32) in matrix representation
are written as follows:
\begin{equation}
\label{e75}
Q \ = \ { \sigma^+ \brack x\sigma^+}\ , \
{\bar Q} \ = \ { x\sigma^- \partial_x-n\sigma^- \brack
-\sigma^-\partial_x} .
\end{equation}
The representation (75) assumes that in the spectral problem (36) an
eigenfunction $\varphi(x)$ is treated  as a two-component spinor
\begin{equation}
\label{e76}
\varphi(x) \ = \ { \varphi_1(x) \brack \varphi_2(x)} ,
\end{equation}
In the matrix formalism, the polynomial space (29) has the form:
\begin{equation}
\label{e77}
{\cal P}_{N,M} \ = \ \left \langle  \begin{array}{c}
x^0,x^1,\dots,x^M  \\
x^0, x^1, \dots,x^N \end{array} \right \rangle
\end{equation}
where the terms of zero degree in $\theta$ come in as the lower component and
the terms of first degree in $\theta$ come in as the upper component.
The operator  ${\bf T}_k (x,\theta)$ becomes a $2 \times 2$ matrix differential
operator ${\bold T}_k (x)$ having
derivatives in $x$ up to $k$-th order. In order to distinguish the matrix
operator in $x$ from the operator in $x,\theta$, we will denote the former as
${\bold T}_k (x)$.
Finally, as a consequence of Theorem 3.1,
 we arrive at the eigenvalue problem for a $2 \times 2$ matrix
quasi-exactly-solvable
differential operator ${\bold T}_k (x)$, possessing in general $(2n+1)$
polynomial solutions of the form ${\cal P}_{n,n-1}$. This eigenvalue
problem can be written in the form (cf. Eq.(37))
\begin{equation}
\label{e78}
 \sum_{i=0}^{i=k} {\bold a}_{k,i} (x) \varphi^{(i)}_{x}
(x) \ = \ \varepsilon \varphi(x) \ ,
\end{equation}
 where the notation $ \varphi^{(i)}_{x}$ means the
$i$th order derivative with respect to $x$ of each component of the spinor
$\varphi(x)$ (see Eq. (76)). The coefficient functions $ {\bold a}_{k,i} (x)$
are given by $2 \times 2$ matrices and generically for the $k$th order
quasi-exactly-solvable operator their matrix elements are polynomials.
Suppose that $k>0$. Then the matrix elements are given
by the following expressions
\begin{equation}
\label{e79}
  {\bold a}_{k,i} (x) \ = \ \left( \begin{array}{cc}
A_{k,i}^{[k+i]} & B_{k,i}^{[k+i-1]} \\
C_{k,i}^{[k+i+1]} & D_{k,i}^{[k+i]}
\end{array}  \right)
\end{equation}
at $k > 0$, where the superscript in square brackets displays the order of
the corresponding polynomial.

 It is easy to calculate the number of free parameters of a
quasi-exactly-solvable operator ${\bold T}_k (x)$
\begin{equation}
\label{e80}
par({\bold T}_k (x))= 4 (k+1)^2
\end{equation}
(cf. Eq.(40)).

For the case of exactly-solvable problems, the matrix elements (79)
of the coefficient functions are modified
\begin{equation}
\label{e81}
  {\bold a}_{k,i} (x) \ = \ \left( \begin{array}{cc}
A_{k,i}^{[i]} & B_{k,i}^{[i-1]} \\
C_{k,i}^{[i+1]} & D_{k,i}^{[i]}
\end{array}  \right)
\end{equation}
where $k > 0$. An infinite family of orthogonal polynomials as
eigenfunctions of Eq. (78), if they exist, will occur, if and only
if the coefficient functions have the form (81). The number of free
parameters of an exactly-solvable operator ${\bold E}_k (x)$ and,
correspondingly, the maximal number of free parameters of the $2 \times 2$
matrix orthogonal polynomials in one real variable, is equal to
\begin{equation}
\label{e82}
par({\bold E}_k (x))= 2k (k+3) + 3
\end{equation}
(cf. Eq.(42)).

The increase in the number of free parameters for the $2 \times 2$ matrix
operators with respect to the case of the operators in $x, \theta$ is
connected to the occurrence of extra monomials of degree $(k+1)$ in
generators of $osp(2,2)$ (see Eqs.(31), (32), (75)), leading to the $k$th
order differential operators in $x$.

 Thus, the above formulas describe the coefficient functions of matrix
differential equations (78), which can possess polynomials in $x$ as
solutions,  resolving the analogue of the
generalized Bochner problem Eq. (0) for the case of $2 \times 2$ matrix
differential equations in one real variable.

Now let us take the quasi-exactly-solvable matrix operator ${\bold T}_2 (x)$
 and try to reduce Eq. (36) to the Schroedinger equation
\begin{equation}
\label{e83}
[ -{1 \over 2} {d^2 \over dy^2} + {\bold V}(y) ] \Psi (y)\ =\ E \Psi (y)
\end{equation}
where $ {\bold V}(y)$ is a two-by-two {\it hermitian} matrix, by making a
change of variable $x \mapsto y$ and \lqq gauge" transformation
\begin{equation}
\label{e84}
\Psi \ = {\bold U} \varphi
\end{equation}
where ${\bold U}$ is an arbitrary $2 \times 2$ matrix depending on
the variable $y$. In order to get some \lqq reasonable" Schroedinger
equation one
should fulfill two requirements: (i) the  potential  ${\bold V}(y)$ must be
hermitian and (ii) the eigenfunctions $\Psi(y)$ must belong to a
certain Hilbert space.

Unlike the case of quasi-exactly-solvable differential operators in one real
variable (see \cite{olver}), this problem has no complete solution so far.
Therefore it seems instructive to display a particular example \cite{st}.

Consider the quasi-exactly-solvable operator
\[ {\bold T}_2 = - 2 T^0 T^- + 2 T^- J - i \beta T^0 Q_1 + \]
\begin{equation}
\label{e85}
\alpha T^0 - (2n+1) T^- - {i\beta \over 2}(3n + 1) Q^1 + {i \over 2}\alpha
\beta Q^2 - i \beta \overline{Q}_1 \ ,
\end{equation}
where $\alpha$ and $\beta$ are parameters. Upon introducing a new
variable $y=x^2$ and  after straightforward calculations
one finds the following expression for the matrix $U$ in Eq. (84)
\begin{equation}
\label{e86}
{\bold U}  = \exp ( - {\alpha y^2 \over 4} + {i \beta y^2 \over 4} \sigma_1)
\end{equation}
and for the potential $\hat{V}$ in Eq. (83)
\[ {\bold V} (y) = {1 \over 8} (\alpha ^2 - \beta ^2) y^2 + \sigma _2
[-(n + {1 \over 4}) \beta + {\alpha \beta \over 4} y^2 - {\alpha \over 4}
\tan {\beta y^2 \over 2}] \cos {\beta y^2 \over 2} +  \]
\begin{equation}
\label{e87}
\sigma _3 [-(n + {1 \over 4}) \beta + {\alpha \beta \over 4} y^2 -
{\alpha \over 4}
\cot {\beta y^2 \over 2}] \sin {\beta y^2 \over 2}
\end{equation}
It is easy to see that the potential ${\bold V}$ is hermitian;
$(2n+1)$ eigenfunctions have the form of polynomials multiplied
by the exponential factor $U$ and they are obviously normalizable.
\section{Polynomials in two real variables}
For further consideration it is convenient to illustrate the space of
polynomials of finite degree in several variables through Newton
diagrams. In order to do this, let us introduce a $d$-dimensional integer
lattice
in ${\bold R}^d$
and put into correspondence to each node with the
coordinates $(k_1,k_2,\ldots k_d)$ the monomial ${x_1}^{k_1}
{x_2}^{k_2} \ldots {x_d}^{k_d}$.
\subsection{Polynomials of the first type.}
Now let us describe the projectivized representation of the
algebra $sl_3({\bold R})$ in the differential operators of the first order
acting on functions of two real variables. It is easy to show that the
generators have the form \cite{st}
\label{e88}
\[ J^1_3\  =\ y^2 \partial_y \  + \ xy\partial_x\  - \ ny \ , \
J^1_2\  =\ x^2 \partial_x \  + \ xy\partial_y\  - \ nx \ ,\]
\[ J^2_3\  =\ -y \partial_x \ ,\ J^2_1\  =\ - \partial_x \ , \
J^3_1\  =\ - \partial_y \ ,\ J^3_2\  =\ -x \partial_y \ , \]
\begin{equation}
 J_d\  =\ y \partial_y \  + \ 2x\partial_x\  - \ n \ , \
\tilde J_d\  =\ 2y \partial_y \  + \ x\partial_x\  - \ n \ ,
\end{equation}
where $x,y$ are the real variables and $n$ is a real number. If $n$ is a
 non-negative integer, the representation becomes finite-dimensional of
dimension $(1+n)(1+n/2)$. The invariant sub-space has a polynomial
basis and is presented as a space of all polynomials of the
following type
\begin{equation}
\label{e89}
{\cal P}^{(I)}_n = \langle 1; x, y; x^2, xy,  y^2;\dots ; x^n, x^{n-1}y,
\dots, x y^{n-1} , y^n \rangle
\end{equation}
or, graphically, the space (89) is given by the Newton diagram of Figure 4.1.
\vskip 1.truecm
\begin{picture}(400,70)(-10,-20)
\noindent\linethickness{1.2pt}
\put(150,10){\line(1,1){30}}
\put(150,10){\line(-1,1){30}}
\linethickness{0.8pt}
\put(120,40){\line(1,0){60}}
\put(152,1){$1$}
\put(110,40){$y^n$}
\put(185,40){$x^n$}
\put(120,40){\circle*{5}}
\put(180,40){\circle*{5}}
\put(150,10){\circle*{5}}
\end{picture}
\begin{center}
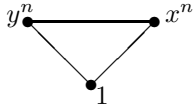
Fig. 4.1 . \ Graphical representation (Newton diagram) of the space
${\cal P}^{(I)}_n$\linebreak (see (89)).
\end{center}
\begin{definition} Let us name a linear differential operator of the $k$-th
order a {\it quasi-exactly-solvable of the first type},${\bf T}_k^{(I)}(x,y)$,
if it preserves the space ${\cal P}^{(I)}_n$ . Correspondingly, the operator
 ${\bf E}_k^{(I)}(x,y) \in {\bf T}_k^{(I)}(x,y)$, which preserves
the infinite flag $ {\cal P}^{(I)}_0 \subset  {\cal P}^{(I)}_1
 \subset {\cal P}^{(I)}_2
\subset \dots \subset {\cal P}^{(I)}_n \subset \dots$ of spaces of all
polynomials of the type (89) , is named an
{\it exactly-solvable of the first type}.
\end{definition}
\begin{LEMMA} {\it Consider the space ${\cal P}^{(I)}_n$.

 (i) Suppose $n > (k-1)$.  Any quasi-exactly-solvable operator of the
first type  ${\bf T}_k^{(I)}(x,y)$, can be
represented by a $k$-th degree polynomial of the generators (88).
 If $n \leq (k-1)$, the part of the quasi-exactly-solvable operator
${\bf T}_k^{(I)}(x,y)$ of the first type containing
derivatives in $x,y$ up to the order $n$ can be represented by a $n$-th
degree polynomial in the generators (88).

(ii) Conversely, any polynomial in (88) is a quasi-exactly solvable operator
of the first type.

(iii) Among quasi-exactly-solvable operators of the first type
there exist exactly-solvable operators of the first type
${\bf E}_k^{(I)}(x,y) \subset {\bf T}_k^{(I)}(x,y)$.}
\end{LEMMA}

{\it Comment 5.} The meaning of the lemma is the following:
${\bf T}_k^{(I)}(x,y)$
at $k < n+1$ is simply an element of the universal enveloping algebra
$U_{sl_3({\bold R})}$ of the algebra $sl_3({\bold R})$ in realization (88).
If $k \geq n+1$, then  ${\bf T}_k (x,y)$ is
represented as a polynomial of $n$th degree in (88) plus
$B {\partial ^{n+1} \over {\partial x^{n-m+1} \partial y^{m}}}$  , where
$m=0,1, \dots (n+1)$ and $B$ is any linear differential operator
of order not higher than $(k-n-1)$.

Let us introduce the grading of the  generators (88) in the following way. The
 generators are characterized by the two-dimensional grading vectors
\label{e90}
\[ deg (J^1_3) = (0,+1) \ , \ deg (J^1_2) = (+1,0) \ , \]
\[ deg (J^2_3) = (-1,+1) \ , \ deg (J^3_2) = (+1,-1) \ , \]
\[ deg (J_d) = (0,0) \ , \ deg (\tilde J_d) = (0,0) \ ,\]
\begin{equation}
deg (J^3_1) = (0,-1) \ , \ deg (J^2_1) = (-1,0) .
\end{equation}
It is apparent that, the grading vector of a monomial in the generators
(88) can be
 defined by the grading vectors of the generators by the rule
\label{e91}
\[ \tilde{deg} (T) \equiv deg [(J^1_3)^{n_{13}} (J^1_2)^{n_{12}}
(J^2_3)^{n_{23}} (J^3_2)^{n_{32}}
(J_d)^{n_{d}} ({\tilde J}_d)^{n_{\tilde d}}(J^3_1)^{n_{31}} (J^2_1)^{n_{21}}]
 \  = \]
\begin{equation}
  (n_{13} + n_{32} - n_{23} - n_{31} \ , \ n_{12} + n_{23} - n_{32} - n_{21})
 \equiv (deg_x (T)\ ,\ deg_y (T))
\end{equation}
Here the $n$'s can be arbitrary  non-negative integers.
\begin{definition} Let us name the grading
of a monomial $T$ in generators (75) the number
\[ deg(T)= deg_x (T) + deg_y (T)\ . \]
We will say that a monomial $T$ possesses positive grading if this number is
positive. If this number is zero, then a monomial has zero grading.
The notion of grading allows one to classify the operators
${\bf T}^{(I)}_k(x,y)$ in a Lie-algebraic sense.
\end{definition}
\begin{LEMMA}{\it A quasi-exactly-solvable operator  ${\bf T}_k^{(I)}
\subset U_{sl_3({\bold R})}$ has no terms of positive grading
 if and only if it is an exactly-solvable operator of the first type.}
\end{LEMMA}

It is worth noting that
among exactly-solvable operators there exists a certain important class of
degenerate operators, which preserve an infinite flag of spaces of all
homogeneous polynomials
\begin{equation}
\label{e92}
\tilde {{\cal P}}^{(I)}_m = \langle  x^m, x^{m-1}y, \dots, x y^{m-1} ,
y^m \rangle
\end{equation}
(represented by a horizontal line in Fig. 4.1).
\begin{LEMMA} {\it A linear differential operator  ${\bf T}_k(x,y)$
preserves
the infinite flag $\tilde {\cal P}^{(I)}_0 \subset \tilde {\cal P}^{(I)}_1
 \subset\tilde {\cal P}^{(I)}_2
\subset \dots \subset\tilde {\cal P}^{(I)}_n \subset \dots$ of  spaces of
all polynomials of the type (89), if and only if it is an
exactly-solvable operator
having terms of zero grading only. Any operators of such a type can be
represented as a polynomial in the generators $J^2_3, J^3_2, J_d,\tilde J_d$
(see (88)), which form the algebra $so_3 \oplus R$ . If such an operator
contains only terms with zero grading vectors, this operator preserves any
space of polynomials.}
\end{LEMMA}
\begin{THEOREM}{\it Let $n$ be a non-negative integer. In general,
the eigenvalue problem for a linear symmetric
 differential operator in two real variables ${\bf T}_k(x,y)$:}
\begin{equation}
\label{e93}
 {\bf T}_k(x,y) \varphi(x,y) \ = \ \varepsilon \varphi(x,y)
\end{equation}
{\it has $(n+1)(n/2+1)$ eigenfunctions in the form of a polynomial in
variables $x, y$ \linebreak belonging to the space (89), if and only if
${\bf T}_k(x,y)$ is a
symmetric quasi-exactly-solvable operator of the first type.
The problem (93) has an infinite sequence of eigenfunctions in the
form of polynomials of the form (89), if and only if the operator
is a symmetric exactly-solvable operator.}
\end{THEOREM}

This theorem gives a general classification of differential equations
\begin{equation}
\label{e94}
\sum_{m=0}^{k}
 \sum_{i=0}^{m} a_{i,m-i}^{(m)} (x,y) {{\partial^m \varphi(x,y)} \over
{\partial x^i \partial y^{m-i} }} \ = \ \varepsilon \varphi(x,y)
\end{equation}
having at least one eigenfunction in the form of a polynomial  in $x,y$ of the
type (89). In general, the coefficient functions $a_{i,m-i}^{(m)} (x,y)$
have quite cumbersome functional structure and we do not display them here
(below we will give their explicit form for ${\bf T}^{(I)}_2(x,y)$). They are
polynomials in $x,y$ of order $(k+m)$ and always contain a general
inhomogeneous polynomial of order $m$ as a part.
The explicit expressions for those polynomials are obtained by substituting
(88) into a general, $k$th order polynomial element of the universal
enveloping algebra $U_{sl_3({\bold R})}$ of the algebra
\linebreak $sl_3({\bold R})$ .
Thus, the coefficients in the polynomials $a_{i,m-i}^{(m)} (x,y)$ can
be expressed through the coefficients of the $k$-th order polynomial element
of the universal enveloping algebra $U_{sl_3({\bold R})}$. The number of free
parameters of the polynomial solutions is defined by the number of parameters
characterizing  a general, $k$th order polynomial element of the universal
enveloping algebra $U_{sl_3({\bold R})}$. In counting free parameters
some relations between generators should be taken into account,
specifically for the given representation (88)
\label{e95}
\[ J^1_2 J_d\ -\ 2 J^1_2 \tilde J_d\ -\ 3 J^1_3 J^3_2\ = \ n J^1_2 \]
\[ J^1_3 \tilde J_d\ -\ 2 J^1_3 J_d\ -\ 3 J^1_2 J^2_3\ = \ n J^1_3 \]
\[ J^3_2 J_d\ +\  J^3_2 \tilde J_d\ -\ 3 J^1_2 J^3_1\ = \ (n+3) J^3_2 \]
\[ J^2_3 J_d\ +\  J^2_3 \tilde J_d\ -\ 3 J^1_3 J^2_1\ = \ (n+3) J^2_3 \]
\[ 3 (J^1_2 J^2_1\  + \ J^1_3 J^3_1\ +\ J^3_2 J^2_3)\ +
  \ J_d J_d\ +  \tilde J_d \tilde J_d \ -\ J_d \tilde J_d\
= \  3 J_d \  +\ 3 n\ +\ n^2 \]
\[ 2J_d J_d\ +\ 2 \tilde J_d \tilde J_d\ -\ 5 J_d \tilde J_d\ +\
9 J^3_2 J^2_3\ =\ (n+6) J_d\ +\ (n-3) \tilde J_d\ +\ n^2\ +\ 3n \]
\[ J_d J_d\ -\  \tilde J_d \tilde J_d\ +\ 3 (J^1_2 J^2_1\  - \ J^1_3 J^3_1)\
 =\ (n+3) (J_d\ -\ \tilde J_d) \]
\[ J_d J^2_1\ -\ 2 \tilde J_d J^2_1\ -\ 3 J^2_3 J^3_1\ = \ n J^2_1 \]
\begin{equation}
 \tilde J_d J^3_1\ -\ 2  J_d J^3_1\ -\ 3 J^3_2 J^2_1\ = \ n J^3_1
\end{equation}
between quadratic expressions in generators (and the ideals generated
by them).\footnote{It was shown that there exist 56 cubic relations
between the generators (88). However, it turned out that all of them
(including standard cubic Casimir operator
in given representation (88)) are
functionally-dependent on quadratic relations (95)
 (Personal communication by G. Post).}
 For the case of exactly-solvable problems, the coefficient functions
$a_{i,m-i}^{(m)} (x,y)$ take the form
\begin{equation}
\label{e96}
a_{i,m-i}^{(m)} (x,y) \ = \ \sum_{p,q=0}^{p+q \leq m} a_{i,m-i,p,q} x^p y^q
\end{equation}
with arbitrary coefficients.

Now let us proceed to the case of the second-order differential equations.
The second-order polynomial in the generators (88) can be represented as
\[ T_2\ =\ c_{\alpha \beta,\gamma \delta} J^{\alpha}_{\beta}
J^{\gamma}_{\delta}
\ +\ c_{\alpha\beta} J^{\alpha}_{\beta}\ +\ c \quad , \]
where we imply summation over all repeating indices; $\alpha, \beta,\gamma,
 \delta$ correspond to the indices of operators in (88) and for the Cartan
generators we suppose both indices simulate $d$ or $\tilde d$, all $c$'s
are set to be real numbers. After substitution the expressions (88) into
$T_2$ the explicit form of the quasi-exactly-solvable operator is given by
\[
{\bf T}_2^{(I)}(x,y)\ =\]
\[ [x^2 P_{2,2}^{xx}(x,y)\ +\ x P_{2,1}^{xx}(x,y)\ +\
\tilde{P}_{2,0}^{xx}(x,y)]{\partial^2 \over \partial x^2}\ +\]
\[ [xy P_{2,2}^{xy}(x,y)\ +\ P_{3,1}^{xy}(x,y)\ +\
\tilde{P}_{2,0}^{xy}(x,y)]{\partial^2 \over \partial x \partial y}\ + \]
\[ [y^2 P_{2,2}^{yy}(x,y)\ +\ y P_{2,1}^{yy}(x,y)\ +\
\tilde{P}_{2,0}^{yy}(x,y)]{\partial^2 \over \partial y^2}\ +\]
\[ [x P_{2,2}^{x}(x,y)\ +\ P_{2,1}^{x}(x,y)\ +\
\tilde{P}_{1,0}^{x}(x,y)]{\partial \over \partial x}\ + \]
\[ [y P_{2,2}^{y}(x,y)\ +\ P_{2,1}^{y}(x,y)\ +\
\tilde{P}_{1,0}^{y}(x,y)]{\partial \over \partial y} +\]
\begin{equation}
\label{e97}
[ P_{2,2}^{0}(x,y)\ +\ P_{1,1}^{0}(x,y)\ +\ \tilde{P}_{0,0}^{0}(x,y)]
\end{equation}
\noindent
where $P_{k,m}^{c}(x,y)$ and $\tilde{P}_{k,m}^{c}(x,y)$ are homogeneous and
inhomogeneous polynomials of order $k$, respectively, the index
$m$ numerates them, the superscript $\underline c$ characterizes the order
of derivative,
which this coefficient function corresponds to. For the case of the
second-order-exactly-solvable operator ${\bf E}_2^{(I)} (x,y)$, the structure
of the coefficient function is similar to (97), except for the fact that all
tilde-less polynomials disappear.
The number of free parameters is equal to
\[par({\bf T}_2^{(I)}(x,y))=36 \ . \]
\noindent
while for the case of an exactly-solvable operator
(an infinite sequence of polynomial eigenfunctions in the problem (93))
\[ par({\bf E}_2^{(I)}(x,y))=25 \ . \]
An important particular case is  when the quasi-exactly-solvable operator
${\bf T}_2^{(I)}(x,y)$ possesses two invariant sub-spaces of the type (89).
This situation is described by the following lemma:
\begin{LEMMA} {\it Suppose ${\bf T}_2^{(I)}(x,y)$ has no terms
of grading 2
\begin{equation}
\label{e98}
c_{12,12}=c_{13,13}=c_{12,13}=0 \ .
\end{equation}
If there exist some coefficient $c$'s and a non-negative integer
 $N$ such that the conditions
\[ c_{12} = (n-N-m)c_{12, d}+(n-2N+m)c_{12,\tilde d}+(N-m)c_{12,32} \ ,\]
\begin{equation}
\label{e99}
 c_{13} = (n-N-m)c_{13, d}+(n-2N+m)c_{13,\tilde d}+ mc_{13,23} \ ,
\end{equation}
are  fulfilled at all $m=0,1,2,\dots,N$, then the operator
${\bf T}_2^{(I)}(x,y)$ preserves both ${\cal P}_n$ and ${\cal P}_N$, and
$par({\bf T}_2^{(I)}(x,y)) = 31$.}
\end{LEMMA}

Now let us proceed to the important item: under what conditions on the
coefficients in (97) the second-order-quasi-exactly-solvable operators can be
reduced to a form of the Schroedinger operator after some gauge
transformation (conjugation)
\begin{equation}
\label{e100}
f(x,y)e^{t(x,y)}{\bf T}_2(x,y)^{(I)} e^{-t(x,y)}\ =\ -\Delta_g\ +\ V(x,y,n)
\end{equation}
where $f, t, V$ are some functions in ${\bold R^2}$, and $\Delta_g$ is the
Laplace-Beltrami operator with some metric tensor $g_{\mu \nu}$;
$ \mu,\nu=1,2$
\footnote{Our further consideration will be restricted the case $f=1$ only.}.
This is a difficult problem which has no complete solution yet.
In \cite{st,gko} a few multi-parametric examples have been constructed.

In general, if one can perform the transformation (100), we arrive at
potentials $V(x,y,n)$ containing explicitly a dependence on the parameter
$n$ and hence the dimension
of the initial invariant subspace (89). Finally, we end up
at quasi-exactly-solvable Schroedinger equations with a certain
number of eigenfunctions known algebraically. Taking an exactly-solvable
operator ${\bf E}_2^{(I)}(x,y)$ instead of the quasi-exactly-solvable one
and performing the transformation
(100), a certain exactly-solvable Schroedinger operator emerges. Evidently,
they will
have no dependence on the parameter $n$ and correspondently an infinite
sequence of eigenfunctions of the Schroedinger equation has a form of
polynomials.

Also it is worth noting that there exists an important particular class of
quasi-exactly-solvable operators, lying in a certain sense in between
quasi-exactly-solvable and exactly-solvable operators.
The algebra $sl_3({\bold R})$ in realization (88) contains the algebra
$so_3({\bold R})$ as a sub-algebra in such a way that
\[ J^1\  =\ (1+ y^2) \partial_y \  + \ xy\partial_x\  - \ ny \ , \quad
J^2\  =\ (1+x^2) \partial_x \  + \ xy\partial_y\  - \ nx \ ,\]
\begin{equation}
\label{e101}
J^3\ =\ x\partial_y \ -\ y\partial_x
\end{equation}
which is not graded anymore.
If the parameter $n$ is a non-negative integer (and coincides with
 that in (88)), the same finite-dimensional invariant sub-space
${\cal P}^{(I)}_{n}$ (see (89) and Fig. 4.1) as in the original
$sl_3({\bold R})$ occurs. For this case a \linebreak
corresponding finite-dimensional representation is reducible
and unitary. It has been proven \cite{st,mprst}, that any symmetric
bilinear combination of generators (101),\
$T_2^{(e)}=c_{\alpha \beta}J^{\alpha}J^{\beta}$,\ where
$c_{\alpha \beta}=c_{\beta \alpha}$ are real numbers,
 can be reduced to a form of the Laplace-Beltrami operator plus a scalar
function making a certain gauge transformation (conjugation).
\footnote{In fact, in \cite{mprst} has been proven a more general
theorem, that any symmetric bilinear combination of generators of unitary
representation of the semisimple Lie algebra, realized in the first-order
differential operators, can be reduced to a form of the Laplace-Beltrami
operator with some metric tensor plus a scalar function making a certain
gauge transformation (conjugation).}
In general, the metric tensor $g_{\mu \nu}$ is not degenerate
\footnote{Degeneracy of $g_{\mu \nu}$ occurs, for example, if
$c_{\alpha \beta}$ has two vanishing
eigenvalues.} and the obtained potential $V(x,y)$ is given
by a rational function of the coordinates. Moreover, the potential $V(x,y)$
has no dependence on the label $n$! All $n$-dependent terms contain no $x,y$-
dependence and can be included into redefinition of the reference point for
eigenvalues. At first sight, those quasi-exactly-solvable operators
look as if they are exactly-solvable, but is not so.
As a consequence of the fact that the
quadratic Casimir operator for $so_3({\bold R})$ in the realization (101)
is non-trivial and commutes with $T_2^{(e)}$, the functional space
of $T_2^{(e)}$ is subdivided into the finite-dimensional blocks
of increasing sizes, corresponding to the irreducible representations
of $so_3({\bold R})$ (unlike truly exactly-solvable operators, for which the
sizes of blocks are fixed and equal to one). We name such operators
{\em exactly-solvable of the second type}.

Another interesting property appears, if the matrix $c_{\alpha \beta}$
has one vanishing eigenvalue: a certain mysterious relation holds \cite{st}
(see also \cite{mprst,lt})
\begin{equation}
\label{e102}
V(x,y,\{c\})\ =\ {3 \over 16} R(x,y,\{c\}),
\end{equation}
where $R(x,y,\{c\})$ is the scalar curvature calculated through the metric
tensor \linebreak attached in the Laplace-Beltrami operator. This can
imply that the
corresponding Schroedinger operator has the form of a purely geometrical
object! The real meaning of this fact is not understood so far.
\subsection{Polynomials of the second type.}
In \cite{t2} we studied the quasi-exactly-solvable operators in one
real variable. It turned out that the solution to this problem was found using
a connection with the projectivized representation of the algebra
$sl_2({\bold R})$.
As a natural step in developing the original idea, let us consider the
projectivized representation of the direct sum of two
algebras $sl_2({\bold R})$.

The algebra $sl_2({\bold R}) \oplus sl_2({\bold R})$ taken in projectivized
representation acts on functions of two real variables.
The generators have the form (see (2))
\label{e103}
\[ J_x^+ = x^2 \partial_x - n x\quad ,\quad J_y^+ = y^2 \partial_y - m y\ , \]
\[ J_x^0 = x \partial_x - {n \over 2} \quad ,\quad
J_y^0 = y \partial_y - {m \over 2}\ , \]
\begin{equation}
J_x^- = \partial_x \quad ,\quad J_y^- = \partial_y \ ,
\end{equation}
where $x,y$ are the real variables and $n,m$ are numbers. If $n,m$ are
non-negative integers, there exists the  finite-dimensional representation
of dimension \linebreak $(n+1)(m+1)$. Evidently, the invariant sub-space
has a polynomial basis and
is presented as a space of all polynomials described by the Newton diagram
of Fig. 4.2. We denote this space as ${\cal P}^{(II)}_{n,m}$.
\vskip 5.truecm
\noindent
\begin{picture}(400,110)(-10,-20)
\linethickness{1.2pt}
\put(150,10){\line(1,1){50}}
\put(150,10){\line(-1,1){30}}
\linethickness{0.8pt}
\put(120,40){\line(1,1){50}}
\put(200,60){\line(-1,1){30}}
\put(152,1){$1$}
\put(105,40){$y^m$}
\put(165,95){$x^ny^m$}
\put(205,60){$x^n$}
\put(120,40){\circle*{5}}
\put(200,60){\circle*{5}}
\put(170,90){\circle*{5}}
\put(150,10){\circle*{5}}
\end{picture}
\begin{center}
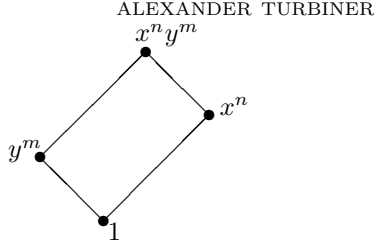
Fig. 4.2. \ Graphical representation (Newton diagram) of the space
${\cal P}^{(II)}_{n,m}$ .
\end{center}
\begin{definition} Let us name a linear differential operator of the $(k,p)$th
order, containing derivatives in $x$ and $y$ up to $k$th
and $p$th orders, respectively\footnote{So a leading derivative has a form
${\partial^{(k+p)} \over {\partial^k_x\partial^p_y}}$. Also we will use a
notation through this section ${\bf T}_{N}(x,y)$ implying that in general all
derivatives of the order $N$ are presented.},
a {\it quasi-exactly-solvable of the second type}, ${\bf T}_{k,p}^{(II)}(x,y)$,
 if it preserves the space ${\cal P}^{(II)}_{n,m}$. \linebreak
Correspondingly,
the operator ${\bf E}_{k,p}^{(II)}(x,y) \in {\bf T}_{k,p}^{(II)}(x,y)$,
which preserves either
the infinite flag $ {\cal P}^{(II)}_{0,m} \subset  {\cal P}^{(II)}_{1,m}
 \subset {\cal P}^{(II)}_{2,m}
\subset \dots \subset {\cal P}^{(I)}_{n,m} \subset \dots$ , or
the infinite flag $ {\cal P}^{(II)}_{n,0} \subset  {\cal P}^{(II)}_{n,1}
 \subset {\cal P}^{(II)}_{n,2}
\subset \dots \subset {\cal P}^{(I)}_{n,m} \subset \dots$
of  spaces of all polynomials,
is named an {\it exactly-solvable of the type} $2_x$ or $2_y$, respectively.
\end{definition}
\begin{LEMMA} {\it Consider the space ${\cal P}^{(II)}_{n,m}$.

 (i) Suppose $n > (k-1)$ and $m > (p-1)$.  Any quasi-exactly-solvable operator
 of the second type ${\bf T}_{k,p}^{(II)}(x,y)$, can be
represented by a $(k,p)$th degree polynomial of the generators (103).
 If $n \leq (k-1)$ and/or $m \leq (p-1)$, the part of the
quasi-exactly-solvable operator ${\bf T}_{k,p}^{(II)}(x,y)$ of the second type
containing derivatives in $x,y$ up to order $n,m$, respectively,
 can be represented by an $(n,m)$th degree polynomial in the generators (103).

(ii) Conversely, any polynomial in (103) is a  quasi-exactly solvable operator
of the second type.

(iii) Among quasi-exactly-solvable operators of the second type
there exist exactly-solvable operators of the second type
${\bf E}_{k,p}^{(II)}(x,y)
\subset {\bf T}_{k,p}^{(II)}(x,y)$.}
\end{LEMMA}

Similarly, as it has been done for the algebra $sl_3({\bold R})$, one can
introduce the notion of grading:
\label{e104}
\[ deg (J_x^+) = (+1,0) \ , \ deg (J_y^+) = (0,+1) \ , \]
\[ deg (J_x^0) = (0,0) \ , \ deg (J_y^0) = (0,0) \ , \]
\begin{equation}
deg (J_x^-) = (-1,0) \ , \ deg (J_y^-) = (0,-1) .
\end{equation}
It is apparent that, the grading vector of a monomial in the generators
(90) can be
 defined by the grading vectors (104) of the generators by the rule
\label{e105}
\[ \tilde{deg} (T) \equiv deg [(J_x^+)^{n_{x+}} (J_x^0)^{n_{x0}}
(J_x^-)^{n_{x-}}
(J_y^+)^{n_{y+}} (J_y^0)^{n_{y0}}(J_y^-)^{n_{y-}}] \  = \]
\begin{equation}
  (n_{x+} - n_{x-} \ , \ n_{y+} - n_{y-}) \equiv (deg_x (T)\ ,\ deg_y (T))
\end{equation}
Here the $n$'s can be arbitrary  non-negative integers.
\begin{definition} Let us name the $x$-grading, $y$-grading and grading
of a monomial $T$ in generators (103) the numbers $deg_x (T)$, $deg_y (T)$ and
 $deg(T)= deg_x (T) + deg_y (T)$, respectively.
We say that a monomial $T$ possesses positive $x$-grading ($y$-grading,
 grading), if the number $deg_x (T)$ ($deg_y (T)$, $deg(T)$) is
positive. If $deg_x (T) (deg_y (T), deg(T)) =0$, then a monomial has
zero $x$-grading ($y$-grading, grading).
\end{definition}

The notion of grading allows one to classify the operators
${\bf T}_{k,p}^{(II)}(x,y)$ in the Lie-algebraic sense.
\begin{LEMMA}{\it The quasi-exactly-solvable operator
${\bf T}_{k,p}^{(II)}(x,y)$ preserves
the infinite flag $ {\cal P}^{(II)}_{0,m} \subset  {\cal P}^{(II)}_{1,m}
 \subset {\cal P}^{(II)}_{2,m}
\subset \dots \subset {\cal P}^{(II)}_{n,m} \subset \dots$
 of spaces of all
the polynomials, if and only if it is an exactly-solvable operator of
type $2_x$ having no terms of positive $x$-grading, $deg_x (T)> 0$.

The quasi-exactly-solvable
operator  ${\bf T}_{k,p}^{(II)}(x,y)$  preserves
the infinite flag $ {\cal P}^{(II)}_{n,0} \subset  {\cal P}^{(II)}_{n,1}
 \subset {\cal P}^{(II)}_{n,2}
\subset \dots \subset {\cal P}^{(II)}_{n,m} \subset \dots$
 of all spaces of
the polynomials, if and only if it is an exactly-solvable operator of
type $2_y$ having no terms of positive $y$-grading, $deg_y (T)> 0$.

 If a quasi-exactly-solvable operator of the second type
contains no terms of
\newpage
\noindent
positive grading, this operator preserves
 the infinite flag $ {\cal P}^{(I)}_0 \subset  {\cal P}^{(I)}_1
 \subset {\cal P}^{(I)}_2
\subset \dots \subset {\cal P}^{(I)}_n \subset \dots$ of spaces of all
the polynomials of the type (76) and is attached to the exactly-solvable
operator of the first type.}
\end{LEMMA}
\begin{THEOREM} {\it Let $n,m$ be non-negative integers. In general,
the eigenvalue problem (93) for a linear symmetric
 differential operator in two real variables ${\bf T}_{k,p}(x,y)$
 has $(n+1)(m+1)$ eigenfunctions in the form of a polynomial in variables
 $x, y$ belonging to the space ${\cal P}^{(II)}_{n,m}$ , if and only if
${\bf T}_{k,p}(x,y)$
is a quasi-exactly-solvable symmetric operator of the second type.
The problem (93) has an infinite sequence of eigenfunctions in the
form of polynomials belonging the space ${\cal P}^{(II)}_{n,m}$ at
fixed $m$ ($n$), if and only if the operator is an exactly-solvable
symmetric operator of type $2_x$ ($2_y$).}
\end{THEOREM}

This theorem gives a general classification of differential equations (94),
having at least one eigenfunction in the form of a polynomial  in $x,y$ of
type ${\cal P}^{(II)}_{n,m}$. In general, the coefficient functions
$a_{i,m-i}^{(m)} (x,y)$ in (94)
have a quite cumbersome functional structure and we do not display them here
 (below we will give their explicit form for ${\bf T}_2^{(II)}(x,y)$).
They are polynomials in $x,y$ of the order $(k+m)$.
The explicit expressions for those polynomials are obtained by substituting
(103) into a general, $k$th order polynomial element of the universal
enveloping algebra $U_{sl_2({\bold R}) \oplus sl_2({\bold R})}$ of the
algebra $sl_2({\bold R}) \oplus sl_2({\bold R})$ .
Thus, the coefficients in the polynomials $a_{i,m-i}^{(m)} (x,y)$ can
be expressed through the coefficients of the $k$-th order polynomial element
of the universal enveloping algebra $U_{sl_2({\bold R}) \oplus
sl_2({\bold R})}$ .
The number of free parameters of the polynomial solutions is defined by the
number of parameters characterizing a general, $k$th order polynomial
element of the universal enveloping algebra
$U_{sl_2({\bold R}) \oplus sl_2({\bold R})}$.
In counting free parameters some relations between generators should be
taken into account, specifically
for the given representation (103)
\label{e106}
\[ J^+_x J^-_x\ -\ J^0_x J^0_x\ +\ J^0_x\
=\ -{n \over 2}({n \over 2}+1) \]
\begin{equation}
 J^+_y J^-_y\ -\ J^0_y J^0_y\ +\  J^0_y\ =\
-{m \over 2}({m \over 2}+1)
\end{equation}
between quadratic expressions in generators (and the ideals generated by them).
\footnote{For this case they correspond to the quadratic Casimir operators.}

Now let us proceed to the case of the second-order differential equations.
The second-order polynomial in the generators (103) can be represented as
\begin{equation}
\label{e107}
T_2\ =\ c_{\alpha \beta}^{xx} J^{\alpha}_{x} J^{\beta}_{x}\ +\
c_{\alpha \beta}^{xy} J^{\alpha}_{x} J^{\beta}_{y}\ +\
c_{\alpha \beta}^{yy} J^{\alpha}_{y} J^{\beta}_{y}\ +\
c_{\alpha}^{x} J^{\alpha}_{x}\ +\
 c_{\alpha}^{y} J^{\alpha}_{y}\ +\ c \quad ,
\end{equation}
where we imply summation over all repeating indices and
$\alpha, \beta = \pm,0$; all $c$'s
are set to be real numbers. Taking (106) in account, it is easy to show that
$T_2$ is characterized by 26 free parameters. Substituting (103) into (107),
we obtain the explicit form of the second-order-quasi-exactly-solvable
operator
\[
{\bf T}_2^{(II)}(x,y)\ =\]
\[  \tilde P_{4,0}^{xx}(x){\partial^2 \over \partial x^2}\ + \
 [x^2y^2 P_{0,2}^{xy}\ +\ xy \tilde P_{1,2}^{xy}(x,y)\ +\
\tilde{P}_{2,0}^{xy}(x,y)]{\partial^2 \over \partial x \partial y}\ + \
 \tilde P_{4,0}^{yy}(y){\partial^2 \over \partial y^2}\ +\]
\begin{equation}
\label{e108}
 [\tilde P_{3,1}^{x}(x)\  +\
y\tilde{P}_{2,0}^{x}(x)]{\partial \over \partial x}\ + \
 [\tilde P_{3,1}^{y}(y)\  +\
x\tilde{P}_{2,0}^{y}(y)]{\partial \over \partial y}\ + \
 \tilde{P}_{2,0}^{0}(x,y)
\end{equation}
\noindent
where $P_{k,m}^{c}(x,y)$ and $\tilde{P}_{k,m}^{c}(x,y)$ are homogeneous and
inhomogeneous polynomials of order $k$, respectively, the index
$m$ marks them, and the superscript $\underline c$ characterizes the order
 of the derivative corresponding to this coefficient function. For the case
of the second-order-exactly-solvable operator of type $2_x$
\[
{\bf E}_2^{(II)} (x,y) \ = \]
\[  \tilde Q_{2,0}^{xx}(x){\partial^2 \over \partial x^2}\ +\
 [x\tilde{Q}_{2,1}^{xy}(y) \ +\ \tilde{Q}_{2,0}^{xy}(y)]
{\partial^2 \over \partial x \partial y}\ + \
 \tilde Q_{4,0}^{yy}(y){\partial^2 \over \partial y^2}\ +\]
\begin{equation}
\label{e109}
 \tilde Q_{1,1}^{x}(x)\  +\
y\tilde{Q}_{1,0}^{x}(x)]{\partial \over \partial x}\ + \
 \tilde Q_{1,1}^{y}(y)\  +\
x\tilde{Q}_{1,0}^{y}(y)]{\partial \over \partial y}\ + \
 \tilde{Q}_{2,0}^{0}(y)
\end{equation}
For the case of the second-order-exactly-solvable operator of type $2_y$
the functional form is similar to (109) with the interchange
$x \leftrightarrow y$.

It is easy to show that the numbers of free parameters are
 \[par({\bf T}_2^{(II)}(x,y))=26 \ . \]
 \[par({\bf E}_2^{(II)}(x,y))=20 \ . \]
As in the previous cases, there is an important particular case of
quasi-exactly-solvable operators of the second order, where they possess two
invariant sub-spaces.
\begin{LEMMA} {\it Suppose in (107) there are no terms of $x$-grading 2
\begin{equation}
\label{e110}
c_{++}^{xx}=0 \ .
\end{equation}
If there exist some coefficient $c$'s and a non-negative integer
 $N$ such that the \linebreak conditions
\[ c_{++}^{xy}=0 \ ,\]
\[ c_{+-}^{xy}=0 \ ,\]
\begin{equation}
\label{e111}
 c_{+}^{x} = (n/2-N)c_{+0}^{xx}+(m/2-k)c_{+0}^{xy} \ ,
\end{equation}
are  fulfilled at all $k=0,1,2,\dots,m$, then the operator
${\bf T}_2^{(II)}(x,y)$
preserves both ${\cal P}_{n,m}^{II}$ and ${\cal P}_{N,m}^{II}$, and
$par({\bf T}_2^{(II)}(x,y)) = 22$.}
\end{LEMMA}

Generically, the question of the reduction of quasi-exactly-solvable
operators of the second type
${\bf T}_2^{(II)}(x,y)$ to the form of the Schroedinger operator is still
open. Recently, several multi-parametrical families of
those Schroedinger operators were constructed in \cite{st,gko}. As well as
the case of the first type of quasi-exactly-solvable operators,
corresponding Schroedinger operators contain
in general the Laplace-Beltrami operator characterizing a non-flat space
metric tensor.

The above analysis of linear differential operators preserving the space
${\cal P}_{n,m}^{(II)}$ can be naturally extended to the case of linear
finite-difference operators defined through the Jackson symbol $D$ (see
Section 2) with the action:
\[ D f(x) = {{f(x) - f(qx)} \over {(1 - q) x}} + f(qx) D \]
instead of a continuous derivative, where $f(x)$ is a real function
and $q$ is a number.
 All above-described results hold
\footnote{ with minor modifications} with replacement of the algebra
$sl_2({\bold R}) \oplus sl_2({\bold R})$ by the quantum algebra
$sl_2({\bold R})_q \oplus sl_2({\bold R})_q$ in \lq projectivized'
representation
(21) (see \cite{t1,t2}). Obviously, in the limit $q \rightarrow 1$ all results
coincide to the results of the present Section.
\subsection{Polynomials of the third type.}
The third case, which we shall discuss here, corresponds to the
family of
Lie algebras  $gl_2 ({\bold R})\ltimes {\bold R}^{r+1}$ (semidirect
sum of  $gl_2 ({\bold R})$ with a $(r+1)$-dimensional abelian ideal; Case 24
in the classification given in \cite{gko1}). This family of
Lie algebras, depending on an integer $r>0$, can be realized in terms of
the first-order differential operators
\label{e112}
\[ J^1\  =\  \partial_x \ , \]
\[ J^2\  =\ x \partial_x\ -\ {n \over 3} \ ,\ J^3\  =\ y \partial_y\ -\ {n
\over {3r}} \ , \]
\[ J^4\  =\ x^2 \partial_x \  +\ rxy \partial_y \ - \ nx \ ,\]
\begin{equation}
 J^{5+i}\  = \ x^{i}\partial_y\ ,\ i=0,1,\dots, r\ ,
\end{equation}
where $x,y$ are real variables and $n$ is a real number.
\footnote{It is worth noting that at $r=1$ the  algebra
$\{ gl_2 ({\bold R})\ltimes {\bold R}^{2}\} \subset sl_3
({\bold R})$. Thus, this case is reduced to one about the first type
polynomials (see Section 4.1). Hereafter, we include $r=1$ into
consideration just for the sake of completeness.}
If $n$ is a non-negative integer, the representation becomes
finite-dimensional. The invariant sub-space has a polynomial basis and
is presented as a space of all polynomials of the form
\begin{equation}
\label{e113}
{\cal P}^{(III)}_{r,n} = \sum_{i,j \geq 0}^{i+rj \leq n} a_{ij} x^i y^j
\end{equation}
or, graphically, (113) is given by the Newton diagram of Figure 4.3. The
general formula for the dimension of the corresponding
finite-dimensional representation (113) is given by
\begin{equation}
\label{e114}
 dim {\cal P}^{(III)}_{r,n}= {{[n^2+(r+2)n+\alpha_{r,n}]} \over {2r}}
\end{equation}
 where for small $r$
\[ \alpha_{1,n}=2 \ , \]
\[ \alpha_{2,n}  = \left\{
\begin{array}{cc}
3 & \mbox{at odd n} \\
4 & \mbox{at even n}
\end{array} \right. \]
\[ \alpha_{3,n}= \left\{
\begin{array}{cc}
4 & \mbox{at (n+1)=0 (mod 3)} \\
6 & \mbox{at other n}
\end{array} \right. \]
\[ \alpha_{4,n}= \left\{
\begin{array}{cc}
5 & \mbox{at (n+1)=0 (mod 4)} \\
8 & \mbox{at other n} \\
9 & \mbox{at (n+3)=0 (mod 4)}
\end{array} \right. \]
\vskip 1.truecm
\begin{picture}(400,120)(-10,-20)
\linethickness{1.2pt}
\put(150,10){\line(1,1){60}}
\put(150,10){\line(-1,1){30}}
\linethickness{0.8pt}
\put(120,40){\line(3,1){90}}
\put(152,1){$1$}
\put(108,41){$y^j$}
\put(215,70){$x^i$}
\put(120,40){\circle*{5}}
\put(210,69){\circle*{5}}
\put(150,10){\circle*{5}}
\end{picture}
\begin{center}
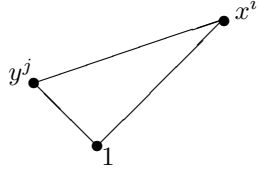
Fig. 4.3 . \ Graphical representation (Newton diagram) of the space
${\cal P}^{(III)}_{r,n}$\linebreak (see (113)).
\end{center}
\begin{definition} Let us name a linear differential operator of the $k$th
order a {\it quasi-exactly-solvable of the $r$-third type},
${\bf T}_k^{(r,III)}(x,y)$, if it preserves the space ${\cal P}^{(III)}_{r,n}$
 . Correspondingly, the
operator ${\bf E}_k^{(r,III)}(x,y) \in {\bf T}_k^{(r,III)}(x,y)$, which
preserves the \linebreak infinite flag $ {\cal P}^{(III)}_{r,0} \subset
 {\cal P}^{(III)}_{r,1} \subset {\cal P}^{(III)}_{r,2}
\subset \dots \subset {\cal P}^{(III)}_{r,n} \subset \dots$ of spaces of all
 polynomials of the type (113) , is named an
{\it exactly-solvable of the $r$-third type}.
\end{definition}
\begin{LEMMA} {\it Take the space ${\cal P}^{(III)}_{r,n}$.

 (i) Suppose $n > (k-1)$.  Any quasi-exactly-solvable operator
${\bf T}_k^{(r,III)}(x,y)$ can be
represented by a $k$th degree polynomial of the generators (112).
 If $n \leq (k-1)$, the part of the quasi-exactly-solvable operator
${\bf T}_k^{(r,III)}(x,y)$ of the $r$-third type containing
derivatives in $x,y$ up to order $n$ can be represented by a $n$-th
degree polynomial in the generators (112).

(ii) Conversely, any polynomial in (112) is a quasi-exactly solvable operator
of the $r$-third type.

(iii) Among quasi-exactly-solvable operators of the $r$-third type
there exist exactly-solvable operators of the $r$-third type
${\bf E}_k^{(r,III)}(x,y) \subset {\bf T}_k^{(r,III)}(x,y)$.}
\end{LEMMA}


One can introduce the grading of the  generators (112) in an analogous way as
has been done before for the cases of the algebra $sl_3({\bold R})$ (see (90))
and $sl_2({\bold R}) \oplus sl_2({\bold R})$ (see (104)). All
 generators are characterized by the two-dimensional grading vectors
\label{e115}
\[ deg (J^1) = (-1,0) \ , \]
\[ deg (J^2) = (0,0) \ , \ deg (J^3) = (0,0) \ , \]
\[ deg (J^4) = (1,0) \ ,\]
\begin{equation}
deg (J^5) = (0,-1) \ , \ deg (J^6) = (1,-1) \ ,\dots, \ deg (J^{5+r}) =
(r,-1)
\end{equation}
Similarly as before, the grading vector of a monomial in the generators
can be defined through the grading vectors of the generators (112)
(cf. (90), (105)).
\begin{definition} Let us name the {\it grading}
of a monomial $T$ in generators (112) the number
\[ deg(T)= deg_x (T) + r\ deg_y (T) \]
 (cf. the case of $sl_3({\bold R})$).
We say that a monomial $T$ possesses positive grading if $deg(T)$ is
positive. If it is zero, then a monomial has zero grading.
\end{definition}

The notion of grading allows us to classify the operators ${\bf T}_k(x,y)$ in
a \linebreak Lie-algebraic sense.
\begin{LEMMA} {\it A quasi-exactly-solvable operator  ${\bf T}_k^{(r,III)}
 \subset U_{gl_2 ({\bf R})\ltimes {\bf R}^{r+1}}$
has no terms of positive grading iff it is an exactly-solvable
operator.}
\end{LEMMA}
\begin{THEOREM} {\it Let $n$ and $(r-1)$ be non-negative integers.
In general, the eigenvalue problem (93) for a linear symmetric
 differential operator in two real variables
${\bf T}_k(x,y)$ has a certain number of eigenfunctions in the form of a
polynomial in variables $x, y$ belonging to the space (113),
if and only if ${\bf T}_k(x,y)$ is a
quasi-exactly-solvable, symmetric operator of the $r$-third type.
The problem (93) has an infinite sequence of eigenfunctions in the
form of polynomials belonging to (113), if and only if the operator is an
exactly-solvable, symmetric operator of the $r$-third type.}
\end{THEOREM}

This theorem gives a general classification of differential equations (94),
having at least one eigenfunction in the form of a polynomial  in $x,y$ of the
type ${\cal P}^{(III)}_{r,n}$. In general, the coefficient functions
$a_{i,m-i}^{(m)} (x,y)$ in (94) are polynomials in $x,y$ and
have a quite cumbersome functional structure and we do not display them here
 (below we will give their explicit form for ${\bf T}_2^{(r,III)}(x,y)$).
The explicit expressions for those polynomials are obtained by substituting
(112) into a general, $k$th order polynomial element of the universal
enveloping algebra $U_{gl_2 ({\bold R})\ltimes {\bold R}^{r+1}}$
of the algebra $gl_2 ({\bold R})\ltimes {\bold R}^{r+1}$.
Thus, the coefficients in the polynomials $a_{i,m-i}^{(m)} (x,y)$ can
be  expressed through the coefficients of the $k$th order
polynomial element of the  universal enveloping algebra
$U_{gl_2 ({\bold R})\ltimes {\bold R}^{r+1}}$.
The number of free parameters of the polynomial solutions is defined by the
number of parameters characterizing a general, $k$th order polynomial element
of the universal enveloping algebra.
In counting free parameters some relations between generators should
be taken into account, specifically
for a given representation (112)
\label{e116}
\[ J^2 J^5\ -\ J^1 J^6\ +\ {n \over 3}J^5\ =\ 0\ , \]
\begin{equation}
 J^1 J^4\ -\ J^2 J^2\ -\ rJ^2 J^3\ -\ J^2\ -\ r({n \over 3}+1) J^3\ =
\ - {n \over 3}({n \over 3}+1) \ ,
\end{equation}
\label{e117}
\[ J^2 J^{6+i}\ +\ rJ^3 J^{6+i}\ -\ J^4 J^{5+i}\ -\ ({n \over 3} +1)J^{6+i}\
=\   0\ ,\]
\begin{equation}
 at\ i=0,1,2,\dots,(r-1) \ ,
\end{equation}
\label{e118}
\[ J^1 J^{7+i}\ -\ J^2 J^{6+i}\ -\ ({n \over 3}+1) J^{6+i}\ =\ 0\ ,\]
\begin{equation}
 at\ i=0,1,2,\dots,(r-2) \ ,
\end{equation}
\label{e119}
\[ J^{5} J^{7+i}\ = \dots =\ J^{5+k} J^{7+i-k}\ , \]
\begin{equation}
 at\ k=0,1,2,\dots \ and\ 2k \leq (2+i)\ , and\ i=0,1,2,\dots ,(2r-2)
\end{equation}
between quadratic expressions in generators (and the ideals generated
by them).

Now let us proceed to the case of the second-order differential equations.
Taking the relations (116)-(119) in account,
one can find the number of free parameters
\label{e120}
\begin{equation}
par({\bf T}_2^{(r,III)} (x,y))\ =\ 5(r+4)
\end{equation}
and obtain the explicit form of the second-order-quasi-exactly-solvable
operator
\[
{\bf T}_2^{(r,III)}(x,y)\ =\]
\[  \tilde P_{4,0}^{xx}(x){\partial^2 \over \partial x^2}\ + \
 [\tilde P_{r+2,1}^{xy}(x)\ +\
y\tilde{P}_{3,0}^{xy}(x)]{\partial^2 \over \partial x \partial y}\ + \
 [\tilde P_{2r,1}^{yy}(x)\ +\
y\tilde{P}_{r+1,0}^{yy}(x)]  ]{\partial^2 \over \partial y^2}\ +\]
\begin{equation}
\label{e121}
 [\tilde P_{3,0}^{x}(x)\ ]{\partial \over \partial x}\ + \
 [\tilde P_{r+1,1}^{y}(x)\  +\
y\tilde{P}_{1,0}^{x}(x)]{\partial \over \partial y}\ + \
 \tilde{P}_{2,0}^{0}(x)
\end{equation}
\noindent
where  $\tilde{P}_{k,m}^{c}(x,y)$ are inhomogeneous polynomials of
order $k$, the index
$m$ numerates them, and the superscript $\underline c$ characterizes the order
 of the derivative corresponding to this coefficient function. For the case
of the second-order-exactly-solvable operator of $r$-third type
\[
{\bf E}_2^{(r,III)}(x,y)\ =\]
\[  \tilde P_{2,0}^{xx}(x){\partial^2 \over \partial x^2}\ + \
 [\tilde P_{r+1,1}^{xy}(x)\ +\
y\tilde{P}_{1,0}^{xy}(x)]{\partial^2 \over \partial x \partial y}\ + \
 [\tilde P_{2r,1}^{yy}(x)\ +\
y\tilde{P}_{r,0}^{yy}(x)]  ]{\partial^2 \over \partial y^2}\ +\]
\begin{equation}
\label{e122}
 [\tilde P_{1,0}^{x}(x)\ ]{\partial \over \partial x}\ + \
 [\tilde P_{r,1}^{y}(x)\  +\
y\tilde{P}_{0,0}^{x}(x)]{\partial \over \partial y}\ + \
 \tilde{P}_{0,0}^{0}(x)
\end{equation}
In this case the number of free parameters is equal to
\label{e123}
\begin{equation}
par({\bf E}_2^{(r,III)} (x,y))\ =\ 5(r+3)
\end{equation}
As in the previous cases, there is an important particular case of
$r$-third type
quasi-exactly-solvable operators, where they possess two invariant sub-spaces.
\begin{LEMMA} {\it Suppose ${\bf E}_2^{(r,III)} (x,y)$ has no
terms of grading 2
\begin{equation}
\label{e124}
c_{44}=0 \ .
\end{equation}
If there exist some coefficient $c$'s and a non-negative integer
 $N$ such that the conditions
\[ c_{4,5+r}=0 \ ,\]
\begin{equation}
\label{e125}
 c_{4}\ =\ ({N \over 3}-m) c_{24}\ +\ ({N \over 3}-k)c_{34}\ ,
\end{equation}
are  fulfilled at all $m,k=0,1,2,\dots$ such that $m+rk=N$, then the
operator ${\bf T}_2^{(r,III)}(x,y)$
preserves both ${\cal P}_{r,n}^{III}$ and ${\cal P}_{r,N}^{III}$, and
$par({\bf T}_2^{(r,III)}(x,y)) =5r+17 $.}
\end{LEMMA}

Generically, the question of the reduction of the quasi-exactly-solvable
operator ${\bf T}_2^{(r,III)}(x,y)$ to the form of the Schroedinger
operator is still open. Initially in \cite{gko} a few
multi-parametrical families of those Schroedinger operators have been
constructed.
Similar to the case of quasi-exactly-solvable operators of the first and the
second types, corresponding Schroedinger operators contain
in general the Laplace-Beltrami operator with non-flat space metric tensor.
\subsection{Discussion.}
Through all the above analysis, one very crucial requirement played a
role: we considered the problem in general position, assuming the
spaces of all polynomials contain polynomials with all possible real
(complex) coefficients. If this requirement is not fulfilled,
then the above connection of the finite-dimensional space of all
polynomials and the finite-dimensional representation of the Lie algebra
is lost and degenerate cases occur. This demands a separate investigation.

So, we described three types of finite-dimensional spaces (see Fig. 4.1-3)
with a basis in polynomials of two real variables, which can be preserved
by linear differential operators. The natural question emerges: is it possible
to find linear differential operators having no representation through
generators of a Lie algebra, which preserve the finite-dimensional space
of polynomials, other than those shown in Fig. 4.1,2,3? The answer
is given by the following conjecture \cite{t5}.
\vskip 0.3truecm
{\bf Conjecture 1.} {\it If a linear differential operator $T_k(x,y)$ does
preserve a finite-dimensional space of all polynomials (presented by
a convex Newton diagram) other than
${\cal P}_{n}^{(I)}, {\cal P}_{n,m}^{(II)}, {\cal P}_{r,n}^{(III)}$, then
this  operator preserves either one of those spaces as well, or
a certain infinite flag of finite-dimensional spaces of all polynomials.}
\vskip 0.3truecm
The general question about reducibility of a second-order polynomial in
generators of the algebras $sl_3({\bold R})$, $sl_2({\bold R})\oplus
sl_2({\bold R})$ and $\{ gl_2 ({\bold R})\ltimes {\bold R}^{r+1}\} ,
r=2,3,4,\ldots$ to the form of the Laplace-Beltrami operator
plus a potential is still open. For all cases where this reducibility
has been performed, a metric tensor appeared which corresponded to
non-flat space. All attempts to find the quasi-exactly-solvable or
exactly-solvable problems in ${\bold R}^2$ corresponding to the flat
space did lead to a problem with separability of variables. We believe
that the following conjecture holds.
\vskip 0.3truecm
{\bf Conjecture 2.} {\it In ${\bold R}^2$ there exist no
quasi-exactly-solvable or exactly-solvable problems containing the
Laplace-Beltrami operator with flat-space metric tensor, which are
characterized by non-separable variables.}
\vskip 0.3truecm

Also the next question, under what conditions the eigenfunctions obtained
belong to some $L_2$-space, is lacking a solution.
The analysis analogous to that which has been performed in
\cite{gko}  for one-dimensional quasi-exactly-solvable problems
should be done.
\section{General considerations}
As the main conclusion one can formulate as the following theorem.
\vskip 0.3truecm
{\bf (Main) Theorem.} {\it Consider a Lie algebra $g$ of the differential
operators of the first order, which possesses a finite-dimensional
irreducible representation ${\cal P}$. Any linear differential operator
acting in ${\cal P}$ (having  ${\cal P}$ as the invariant subspace) can
be represented by a polynomial in generators of the algebra $g$ plus an
operator annihilating  ${\cal P}$.}
\vskip 0.3truecm

As a probable extension of this theorem, we assume that the following
conjecture holds \cite{t5}.
\newpage
{\bf Conjecture 3.} {\it If a linear differential operator $T$ acting
on functions in ${\bold R^k}$ does possess a single finite-dimensional
invariant subspace with polynomial basis (presented by a convex Newton
diagram), this finite-dimensional space coincides with a certain
finite-dimensional representation of some Lie algebra, realized
by differential operators of the first order.}
\vskip 0.3truecm
The interesting question is how to describe finite modules of smooth
functions in ${\bf R^k}$ which can serve as invariant sub-spaces of linear
operators.
Evidently, if we could find some realization of a Lie algebra $g$ of
differential operators, possessing an irreducible finite module of smooth
functions, we can immediately consider the direct sum of several species of
$g$ acting on the functions given at the corresponding direct products of
spaces, where the original algebra $g$ acts in each of them. As an example,
this procedure is presented in Section 4.2: the direct sum of two
$sl_2({\bold R})$ generates the linear differential operators acting on the
functions at ${\bold R^2}$ having
a rectangular ${\cal P}^{(II)}_{n,m}$ as the invariant sub-space. Taking
the direct sum of $k$ species of $sl_2({\bold R})$ acting at ${\bold R^k}$,
we arrive to the operators possessing a
$k$-dimensional parallelopiped as the invariant sub-space.
In this case the algebra $sl_2({\bold R})$ plays the role of a
\lqq primary" algebra giving rise to a multidimensional convex geometrical
 figure (Newton diagram) as the invariant
sub-space of some linear differential operators. In the space
${\bold R^1}$ there is only one such  \lqq irreducible" convex Newton
diagram - a finite interval, in other words, a space of all polynomials
in $x$ of degree not higher than a certain positive integer.
Thus, in ${\bold R^2}$ the rectangular Newton diagram is \lqq reducible"
stemming from the direct product of two intervals as one-dimensional
Newton diagrams.
In turn, the irreducible convex Newton diagrams in ${\bold R^2}$ are
likely exhausted by different triangles, connected to the algebras
$sl_3({\bold R})$ and $ gl_2 ({\bold R})\ltimes {\bold R}^{r+1} ,
r=2,3,4,\ldots$ (see Figures 4.1, 4.3). Those algebras
play a role of primary algebras. At least,
we could not find any other convex Newton diagrams in ${\bold R^2}$
(see conjecture 1).

Before description of our present knowledge about convex Newton
diagrams
for the ${\bold R^k}$ case, let us recall the regular representation
of the Lie algebra $sl_N({\bold R})$ given on the flag manifold
which acts on the smooth functions of $N(N-1)/2$ real \linebreak
variables $z_{i,i+q}\ , i=1,2,\dots,(N-1),\ q=1,2,\dots,(N-i)$.
The explicit formulas for generators in the Chevalley basis are given by
(see \cite{r})
\label{e126}
\[ D(e_i)=-{\partial \over {\partial z_{i,i+1}}}- \sum_{q=i+2}^{N}
z_{i+1,q}{\partial \over {\partial z_{i,q}}} \]
\[ D(f_i)=\sum_{q=1}^{i}z_{q,i+1}z_{i,i+1}{\partial \over {\partial z_{q,i+1}}}
+
\sum_{q=1}^{i-1}(z_{q,i+1}-z_{q,i}z_{i,i+1}){\partial \over {\partial z_{q,i}}}
\]
\[
- \sum_{q=i+2}^{N}z_{i,q}{\partial \over {\partial z_{i+1,q}}} \ - \
n_i z_{i,i+1} \]
 \[ D( h_i)=-\sum_{q=i+1}^{N} z_{i,q}{\partial \over {\partial z_{i,q}}}+
\sum_{q=1}^{i-1} z_{q,i}{\partial \over {\partial z_{q,i}}}+ \]
\begin{equation}
\sum_{q=i+2}^{N} z_{i+1,q}{\partial \over {\partial z_{i+1,q}}}-
\sum_{q=1}^{i} z_{q,i+1}{\partial \over {\partial z_{q,i+1}}}
+ n_{i}
\end{equation}
\noindent
where we use a notation $e_i,f_i, h_i,\ i=1,2,\dots,(N-1)$ for generators
of positive and negative roots, and the Cartan generators, respectively.
The algebra $sl_N({\bold R})$ is generated by generators
$e_i,f_i, h_i$. If $n_i$ are
non-negative integers, the finite-dimensional irreducible representation of
$sl_N({\bold R})$ will occur in the form of inhomogeneous polynomials in
variables $z_{i,i+q}$. The highest weight vector is characterized by the
integer numbers $n_i,\ i=1,2,\dots, (N-1)$. If some integers $n_i$ become
zero, the representation (126) is degenerated and acts on a space of
dimension lower than $N(N-1)/2$ (see, e.g., (88) as an example of a
degenerate representation of $sl_3({\bold R})$).

In connection with the space ${\bold R^3}$, at least two irreducible convex
Newton \linebreak diagrams appear:
a tetrahedron, related to a degenerate finite-dimensional \linebreak
representation
of $sl_4({\bold R})$ (the highest weight vector is characterized by two
zeroing and one non-zeroing components), and certain geometrical bodies,
corresponding to the regular representation of $sl_3({\bold R})$ on
 the flag manifold (see (126) at $N=3$) and a representation of
$sp_2({\bold R})$, respectively. We do not know whether those
exhaust all irreducible Newton diagrams or not.

For the general case ${\bold R^k}$, firstly, we wish to indicate, that
there also exist Newton \linebreak polyhedra related to the almost
degenerate representation of
$sl_{k+1}({\bold R})$  (in formulas (126)  all integers except the one
characterizing the highest weight vector become zero). Also
there exist some convex geometrical figures
related to some finite-dimensional degenerate and, sometimes, regular
representations of $sl_{k+1-i}({\bold R})$,\linebreak $i>0$.\footnote{In
order to decide what algebras act on ${\bold R^k}$, it is
necessary to perform some analysis of the Young tableau. I am grateful
to W. R\"uhl for an enlighting discussion of this question.}
It is easy to show that the algebra $so_{k+1}({\bold R})$ acts on
${\bold R^k}$ in a degenerate representation
\label{e127}
\[ H_{ij}\ =\ x_i {\partial \over {\partial x_j}} \ -\
x_j {\partial \over {\partial x_i}} \ , \ i>j \
, \ i,j=1,2,\ldots,k \]
\begin{equation}
G_i\ =\ \sum_{j=1,i \neq j}^{k} x_i x_j {\partial \over {\partial x_j}}
\ + \ (1+x_i^2){\partial \over {\partial x_i}} \ -\ nx_i \
,\ i=1,2,\ldots,k
\end{equation}
where $x_i, i=1,2,\ldots,k$ are real variables and $n$ is a number.
If $n$ is a non-negative integer, the algebra (127) possesses the
finite-dimensional invariant sub-space,
\footnote{The corresponding finite-dimensional representation is unitary
and reducible (see, e.g., the discussion in the end of Section 4.1 and
also in \cite{mprst}).}  which has the form of a
polyhedron in ${\bold R^k}$ and coincides with a finite-dimensional
invariant sub-space of the algebra $sl_{k+1}({\bold R})$ in the almost
degenerate representation, acting on ${\bold R^k}$.
Making a gauge transformation, any symmetric bilinear combination in the
generators (127) of the algebra $so_{k+1}({\bf R})$ with real coefficients
can be reduced to the form of a Laplace-Beltrami operator with some metric
tensor plus a scalar function giving rise to exactly-solvable
operators of the second type (for the proof see
\cite{mprst} and the discussion in the end of Section 4.1).

The above procedure can be extended to the case of the Lie super-algebras
and also the quantum algebras. For the former case, our experience is
very limited. The only statement one can make is that for linear operators
acting on functions of one real and one Grassmann variable (or,
equivalently, on two-component spinors in on real variable)
there can appear finite-dimensional invariant sub-spaces other than
finite-dimensional representations  of the algebra $osp(2,2)$
(see Section 3.1 and the discussion in \cite{t3}). Nevertheless, those
sub-spaces are connected with the
representations of polynomial elements of the universal enveloping algebra
of the algebra $osp(2,2)$
and, correspondingly, they can be constructed through finite-dimensional
representations of the original super-algebra $osp(2,2)$. Probably, there
are no more irreducible convex Newton diagrams in this space.
For the latter (quantum algebras), there exists an irreducible Newton
diagram in ${\bold R}$
 -- an interval (see Section 2) -- but we could not find any irreducible
(non-trivial) Newton diagrams in ${\bold R^k}$ for $k>1$; however,
at least, in ${\bold R^2}$ they occur, if a quantum plane is introduced,
$xy=qyx$, instead of an ordinary one, $xy=yx$, or a quantum hyperplane
\cite{wz} for the case of more than two variables
\footnote{I am grateful to O.V. Ogievetsky for an interesting
discussion of this subject and valuable comment.}.

It is worth noting that
one can construct quasi-exactly-solvable and exactly-solvable
operators considering \lqq mixed" algebras: a direct sum of the Lie
algebras, Lie superalgebras and the quantum algebras, realized
in the first-order differential and/or difference operators,
possessing finite-dimensional invariant subspaces. In this case the
Newton diagrams are reducible.

As for the cases considered in Sections 1-4, the problem of
reducibility of \linebreak general quasi-exactly-solvable and
exactly-solvable operators to the Schroedinger-type operators is open.

Above we have described the linear operators which possess a
finite-dimensional invariant sub-space with a polynomial basis.
Certainly, making changes of \linebreak variables and
gauge transformations one can obtain the operators with an invariant
sub-space with a non-polynomial basis (but emerging from a polynomial
one, see e.g. (11)). Those  `induced' linear spaces we consider as the
equivalent to the original, polynomial ones.
A general problem can be posed: is it possible to
describe all linear differential
operators possessing a finite-dimensional invariant sub-space
with a basis in a certain explicit form\footnote{It is quite difficult
to define precisely the meaning of the words \lqq space
with a basis in a certain explicit form." One of the examples is the
linear space $<1, x, f_1(x), f_2(x),\ldots , f_n(x)>$,
where $f_i, i=1,2,\ldots n$ are any functions on $R^1$. Other examples
of such spaces are the spaces (1), (11), (29), (77), (89), (113).}, which
 {\it cannot} be reduced to a basis of polynomials using a
change of variables and a gauge transformation?
Equivalently, this problem can be re-formulated as a problem of a
description of finite-dimensional linear spaces of functions, where some
linear differential operators can act. Since the original aim of
the present investigation was mainly the
construction of the quasi-exactly-solvable Schroedinger operators
\cite{t4,st,mprst,olver,gko,lt}, we did not touch on the general problem
(except for a short remark on page 4).
A solution of this problem plays a crucial role for a complete
classification of the  Schroedinger operators
possessing a finite-dimensional invariant sub-space with an explicit
basis in functions. Although it is rather surprising, it is very likely
that this general and abstract problem has a constructive solution, at least
for the one-dimensional Schroedinger operators. This work is in progress.
\footnote{Recently, a classification of linear differential operators
possessing a finite-dimensional invariant sub-space with a basis in monomials
 was completed \cite{pt}.}

\noindent
{\bf Concluding remarks.}

(I).\quad
A crucial approach to the problem of classification of linear operators
possessing a finite-dimensional invariant sub-space(s) with a basis in
functions in an explicit form\footnote{see footnote 18} is connected to a
theory of Riemann surfaces.

Take a symmetric operator depending on a
real parameter polynomially. and possessing the infinite discrete set
of eigenvalues. Searching of eigenvalues is equivalent to solving the
transcendental secular equation. It is rather natural that the roots of
this equation (eigenvalues) form an infinite-sheet Riemann surface as a
parameter manifold. Generically, this surface is characterized by
square-root branch points (see the excellent paper by Bender and Wu
\cite{bw} and also \cite{si,sh,tu}). Existance of $N$-dimensional invariant
sub-space with an
explicit basis implies, that the secular equation is reduced to two
equations: an algebraic one and a transcendental one. It reflects the fact,
that the original determinant factors into the product of two determinants.
Equivalently, the infinite-sheet Riemann surface is split off into two
disjoint Riemann
surfaces: $N$-sheet one and an infinite-sheet one. If there exist $k$
finite-dimensional invariant sub-spaces (for instance, which are
embedded consequently one into the other, or, another case, any two
have no intersection), the original,
infinite-sheet Riemann surface is split off into $k$ finite-sheet Riemann
surfaces and an infinite-sheet one. Namely, such a situation occurs for a
general quasi-exactly-solvable operator. As for exactly-solvable operators,
the corresponding infinite-sheet Riemann surface is split off completely
for separate sheets\footnote{Of course, if the eigenfunctions related with
the infinite flag of invariant sub-spaces form a complete set of the
eigenfunctions of an original operator. In contrary, an additional Riemann
surface is left as well.}. Following this idea, in \cite{t4,tlame}
a complete classification was done of second-order differential
operators in ${\bold R}$ related with $sl_2({\bold R})$-hidden algebra.
This insight allowed to separate the quasi-exactly-solvable and
exactly-solvable operators from generic operators, possessing abstract
finite-dimansional invariant sub-spaces. I hope that such an approach
can lead to an abstract classification of linear operators.

(II).\quad The above-discussed questions about the
existence of convex Newton diagrams (equivalently, the
finite-dimensional linear spaces with polynomial basis) as an invariant
sub-space of some linear differential
operator are related to a classification of finite-dimensional Lie
algebras of first-order differential operators possessing an invariant
sub-space. Originally, this problem was formulated by Sophus Lie, who
solved this problem for the algebras on ${\bold C^1}$ and gave a complete
classification of the algebras of vector fields on ${\bold C^2}$ (see
\cite{gko2} and references therein). A complete classification of Lie algebras
of first-order differential operators on ${\bold C^2}$ has been presented
just recently \cite{gko1}, while the same problem for the algebras on
${\bold R^2}$ has not been solved so far, although the vector fields on
 ${\bold R^2}$ have been described \cite{gko2}. The same problem for algebras
on ${\bold C^3}$ has also been discussed by S. Lie, who presented 35 algebras
of vector fields, emphasizing that this list is not complete. As far as we
know at present, similar questions about a classification of the
Lie super-algebras, both vector fields and first-order differential
operators in real and Grassmann variables, are never discussed in literature.
The same situation takes place for the algebras of first-order
finite-difference operators on the real line, or a quantum hyperplane
\cite{wz} .
\vskip 1.2truecm
\begin{center}
{\bf ACKNOWLEDGEMENT}
\end{center}
\vskip 0.4truecm
In closing, I am very indebted to V. Arnold, J. Froehlich, L. Michel,
G. Post, W. R\"uhl and M. Shubin for their interest in the subject and
numerous valuable discussions. Also I am very grateful to the Institute
of Theoretical Physics, ETH-Honggerberg, Zurich, where this work
was mainly done, and the Research Institute for Theoretical
Physics, University of Helsinki and the Institut des Hautes \'Etudes
Scientifiques, Bures-sur-Yvette, where this work was completed,
for their kind hospitality extended to me. I am deeply appreciate to
N. Kamran and P. Olver for their kind invitation to make a contribution
to present volume.

\renewcommand{\theequation}{A.{\arabic{equation}}}
\setcounter{equation}{0}
\newpage
\begin{center}
{\bf APPENDIX}
\end{center}
\vskip 0.4truecm

In connection to the finite-dimensional representations of the algebras
of first-order differential (difference) operators some operator identities
occur \cite{tp}.

1. The following {\it operator identity} holds:
\label{ae1}
\begin{equation}
{(J^+_n)}^{n+1} \equiv (x^2 \partial_x - n x)^{n+1} =
x^{2n+2}\partial_x^{n+1}\ , \ \partial_x \equiv {d \over dx},
\ n=0,1,2,\ldots
\end{equation}
\begin{pf}The proof is straightforward:
\begin{enumerate}
\item[(i)] the operator ${(J^+_n)}^{n+1}$
annihilates the space of all polynomials of power not higher than $n$,
${\cal P}_n(x)=Span\{x^i: 0 \leq i \leq n\}$;
\item[(ii)] in general, an $(n+1)$-th order linear differential
operator annihilating ${\cal P}_n(x)$ must have the form
$B(x)\partial_x^{n+1}$, where $B(x)$ is an arbitrary function and
\item[(iii)]
since ${(J^+_n)}^{n+1}$ is a graded operator, deg$(J^+_n)=+1$,
\footnote{so $J^+_n$ maps $x^k$ to a multiple of $x^{k+1}$ }
deg$({(J^+_n)}^{n+1})=n+1$, hence $B(x)=b x^{2n+2}$, while clearly the
constant $b=1$.
\end{enumerate}
\end{pf}
It is worth noting that taking the power in (A.1) different from
 $(n+1)$, the l.h.s. in (A.1), generally, will contain all derivative
terms from zero up to $(n+1)$-th order with non-vanishing coefficients.

The identity (A.1) has a Lie-algebraic interpretation.
The operator $(J^+_n)$ is the positive-root generator of the algebra $sl_2$
of first-order differential operators (the other $sl_2$-generators are
$J^0_n = x \partial_x - n/2 \ , J^-_n = \partial_x $, see (2)).
Correspondingly, the
space ${\cal P}_n(x)$ is nothing but the $(n+1)$-dimensional
irreducible representation of $sl_2$. The identity (A.1) is a
consequence of the fact that ${(J^+_n)}^{n+1}=0$ in matrix
representation \footnote{It is a particular case of more general
statement: {\it in the algebra $gl_k$
a positive (negative) -root generator $J_{\pm\alpha}$ in a
finite-dimensional representation of the dimension $d$, taken in power  of
of the dimension $d$ is equal to zero,} $(J_{\pm\alpha})^d=0$.}.

Another Lie-algebraic interpretation of (A.1) is connected with occurence of
some relations between elements of the universal enveloping algebra of
three-dimensional Heisenberg algebra, $H_1$: $\{ P, Q, 1\}$.
Let $[P,Q]=1$, then
\label{ae2}
\begin{equation}
(Q^2 P - n Q)^{n+1} = Q^{2n+2} P^{n+1}, \ n=0,1,2,\ldots
\end{equation}
and certainly this identity is more general than (A.1). Taking different
representations for $P$ and $Q$ in terms of differential operators
(others than the standard one, $P=\partial_x$ and $Q=x$),
one can get various families of operator identities, other than (A.1).
Also the formula (A.2) has a meaning of a formula of an ordering:
$Q$-operators are placed on the left, $P$-operators on the right.

One can easily check, that once two operators obey $[P,Q]=1$, then
\[ J^+\ =\ Q^2 P - n Q \ ,\]
\begin{equation}
 J^0\ =\ Q P - n/2 \ ,
\end{equation}
\[ J^-\ =\ P \ , \]
obey $sl_2$-algebra commutation relations. The representation (2) is a
particular case of this representation.
In fact, (A.3) is one of possible embeddings of the algebra $sl_2$ into
the universal enveloping algebra of three-dimensional Heisenberg algebra.

There exist other algebras of differential or finite-difference
operators (in more than one variable), which admit a
finite-dimensional representation.  This leads to more general
and remarkable operator identities. In fact, (A.1) is one representative
of an infinite family of identities for differential and finite-difference
operators. Also there occur generalizations of the relation (A.2) for
a certain polynomial elements of the universal enveloping algebra
of $(2n+1)$-dimensional Heisenberg algebra.

2. The Lie-algebraic interpretation presented above allows us to generalize
(A.1) for the case of differential operators of several variables,
taking appropriate powers of the highest-positive-root generators of
(super) Lie algebras of first-order differential operators, possessing
a finite-dimensional invariant sub-space.
First we consider
the case of $sl_3$. There exists a representation of $sl_3(\bold R)$ as
differential operators on $\bold R^2$ (see Section 4.1, eq.(88)).
One of the generators is
\[ J^1_2 (n)= x^2 \partial_x\ +\ xy \partial_y - n x  \ .\]
The space (89), ${\cal P}_n(x,y)=Span\{x^iy^j: 0 \leq i+j \leq n\}$ is a
finite-dimensional representation for $sl_3$, and due to the fact that
$(J^1_2 (n))^{n+1}=0$ on the space ${\cal P}_n(x,y)$, we have
\label{ae3}
\begin{eqnarray}
{(J^1_2 (n))}^{n+1} \equiv (x^2 \partial_x\ +\ xy \partial_y - n x)^{n+1} =
\nonumber \\
\sum_{k=0}^{n+1} {n+1 \choose k} x^{2n+2-k}y^k \partial_x^{n+1-k}
\partial_y^k  \ .
\end{eqnarray}
This identity is valid for $y \in \bold R$ (as described above), but also if
$y$ is a Grassmann variable, i.e. , $y^2=0$ \footnote{In this case just
two terms in the l.h.s. of (A.4) survive.}. In the last case, $J^1_2 (n)$ is
a generator $T^+$ of the algebra $osp(2,2)$ (see (31)).

More generally (using algebra $sl_{k+1}$ instead of $sl_3$), the following
operator identity holds:
\label{ae4}
\begin{eqnarray}
{(J^{k-1}_k (n))}^{n+1} \equiv (x_1(\sum_{m=1}^{k} x_m \partial_{x_m}\
- n))^{n+1} =
\nonumber \\
x_1^{n+1}\sum_{j_1+j_2+\ldots+j_k=n+1} C^{n+1}_{j_1,j_2,\ldots,j_k}
x_1^{j_1}x_2^{j_2}\ldots x_k^{j_k} \partial_{x_1}^{j_1}
\partial_{x_2}^{j_2}\ldots \partial_{x_k}^{j_k}   \ ,
\end{eqnarray}
where $C^{n+1}_{j_1,j_2,\ldots,j_k}$ are the generalized binomial
(multinomial) coefficients.
If $x \in \bold R^k$, then $J^{k-1}_k (n)$ is a generator of the algebra
$sl_{k+1}(\bold R)$, generated by $D(f_i)$ (see (126)).
While some of the variables $x$'s are Grassmann
ones, the operator $J^{k-1}_k (n)$ is a generator of a certain super
Lie algebra of first-order differential operators. The operator in the
l.h.s. of (A.5) annihilates the linear space of polynomials
${\cal P}_n(x_1,x_2,\ldots x_k)=Span\{x_1^{j_1}x_2^{j_2}\ldots x_k^{j_k} :
 0 \leq j_1+j_2+\ldots+j_k \leq n\}$,
which is represented by tetrahedron Newton diagram in $\bold R^k$.

Considering $(2k+1)$-dimensional Heisenberg algebra:
\[
[P_i,Q_j]=\delta _{ij},\ [P_i,P_j]=0, \ [Q_i,Q_j]=0,\ i,j=1,2,\ldots k
\]
one can rewrite (A.5) in more general form
\label{ae5}
\begin{eqnarray}
(Q_1(\sum_{m=1}^{k} Q_m P_m \ -\ n))^{n+1} =
\nonumber \\
Q_1^{n+1} \sum_{j_1+j_2+\ldots+j_k=n+1} C^{n+1}_{j_1,j_2,\ldots,j_k}
Q_1^{j_1}Q_2^{j_2}\ldots Q_k^{j_k} P_1^{j_1}P_2^{j_2}\ldots P_k^{j_k}   \ ,
\end{eqnarray}
(cf.(A.2)). It is clear, the formulas for $sl_{k+1}$-generators at $k>1$
(see (126)) in terms of $P_i, Q_i$, analogous to (A.3), can be easily derived.

3. For the case of the algebra $ gl_2 ({\bold R})\ltimes {\bold R}^{r+1} ,
r=1,2,\ldots$ (see (112)) the situation \linebreak becomes rather
sophisticated
\begin{eqnarray}
\label{ae6}
(J^4_n)^{n+1} \equiv (x^2 \partial_x \  +\ rxy \partial_y \ - \ nx)^{n+1} =
\nonumber \\
\sum_{k=0}^{n+1} r^k {n+1 \choose k} x^{2n+2-k} y^{n+1-k}\partial_x^k
\partial_y^{n+1-k} + Q_{r,n} \ ,
\end{eqnarray}
where $Q_{1,n}=0$ and also
\[ Q_{r,n} \ =\ \left\{
\begin{array}{cc}
0 & \mbox{at n=0 } \\
r(r-1)x^2 y\partial_y & \mbox{at n=1 } \\
r(r-1)x^3y[3ry \partial_y^2 + 3x \partial_{xy}^2 + (r-2) \partial_y] &
\mbox{at n=2 }
\end{array} \right. \]
(cf. (A.4)). For arbitrary $n$ and $r$
\[ Q_{r,n} = {{r(r-1)n(n+1)} \over 2} \]
\[ [x^{2n}y \partial_x^{n-1} \partial_y +
r(n-1)x^{2n-1}y^2 \partial_x^{n-2} \partial_y^2+
{r^2 \over 4} (n-1)(n-2)x^{2n-2}y^3 \partial_x^{n-3} \partial_y^3+\ldots]
\]
\[ + \ [lower\ order \ derivatives] . \]

The r.h.s. of (A.7) still has quite regular structure,
which can be analyzed using the following trick: variables $y$'s in the
l.h.s. of (A.4), (A.7) can be related to each other through the very
simple connection
\[ y_{(A.7)} \ =\ y_{(A.4)}^r \ ,\]
and then $J^1_2$ (see (A.4)) coincides with $J^4_n$.

4. The above-described family of operator identities can be generalized
for the case of finite-difference operators with the Jackson symbol, $D_x$
(see, e.g., \cite{e}), obeyed the following Leibnitz rule
\[ D_x f(x) = {{f(x) - f(q^2x)} \over {(1 - q^2) x}} + f(q^2x) D_x \ ,\]
instead of the ordinary, continuous derivative.\footnote
{The definition of the Jackson symbol used here is slightly different
than that presented in Section 2.} Here, $q$ is an arbitrary
complex number. The following operator identity holds:
\label{ae7}
\begin{equation}
{(\tilde J^+_n)}^{n+1} \equiv  (x^2 D_x - \{ n \}  x)^{n+1} =
q^{2n(n+1)}x^{2n+2} D^{n+1}_x , n=0,1,2,\ldots
\end{equation}
(cf. (A.1)), where $\{n\} = {{1 - q^{2n}}\over {1 - q^2}}$ is a so-called
$q$-number.  The operator in the l.h.s. annihilates the space (1), ${\cal
P}_n(x)$. The proof is similar to the proof of identity (A.1).

From the algebraic point of view the operator $\tilde J^+_n$ is the generator
of a $q$-deformed algebra $sl_2(\bold R)_q$ of first-order
finite-difference operators on the line: \linebreak
$\tilde J^0_n = \ x D - \hat{n},\ \tilde J^-_n = \ D $, where
$\hat n \equiv {\{n\}\{n+1\}\over \{2n+2\}}$ (see \cite{ot} and
also \cite{t1}), obeying the commutation relations (22).

One can show, that once two operators obey
$\tilde P \tilde Q-q^2 \tilde Q \tilde P=1$, then
\begin{equation}
(\tilde Q^2 \tilde P - \{n\}\tilde Q)^{n+1} = q^{2n(n+1)}\tilde Q^{2n+2}
\tilde P^{n+1}, \ n=0,1,2,\ldots
\end{equation}
holds and
\[ \tilde J^+\ =\ \tilde Q^2 \tilde P -\{ n\} \tilde Q \ ,\]
\begin{equation}
\tilde J^0\ =\ \tilde Q \tilde P - \hat{n} \ ,
\end{equation}
\[ \tilde J^-\ =\ \tilde P \ , \]
obey $sl_2(\bold R)_q$-algebra commutation relations (22).
It is evident, that this representation is more general than the
representation (21).

An attempt to generalize (A.4) replacing continuous derivatives by Jackson
symbols immediately leads to requirement to introduce the quantum plane and
$q$-differential calculus \cite{wz}
\label{ae8}
\[ xy=qyx\ , \]
\[ D_x x=1+q^2 xD_x+(q^2-1) yD_y \quad ,\quad D_x y=qyD_x \ ,\]
\[ D_y x=qxD_y\quad , \quad D_y y=1+q^2 yD_y \ ,\]
\begin{equation}
 D_xD_y=q^{-1}D_yD_x \ .
\end{equation}
The formulae analogous to (A.4) have the form
\label{ae9}
\[ {(\tilde J^1_2 (n))}^{n+1} \equiv
(x^2 D_x\ +\ xy D_y - \{ n\} x)^{n+1} = \]
\begin{equation}
\sum_{k=0}^{n+1} q^{2n^2-n(k-2)+k(k-1)} {n+1 \choose k}_q x^{2n+2-k}y^k
D_x^{n+1-k} D_y^k  \ ,
\end{equation}
where
\[ {n \choose k}_q \equiv {\{n\}! \over {\{k\}!\{n-k\}!}}\ ,\ \{n\}! =
\{1\} \{2\}\ldots \{n\} \]
are $q$-binomial coefficient and $q$-factorial, respectively. As in all
previous cases, if $y \in {\bold R}$, the operator $\tilde J^1_2 (n)$ is
one of generators of the
$q$-deformed algebra $sl_3(\bold R)_q$ of finite-difference
operators, acting on the quantum plane and having the linear space
${\cal P}_n(x,y)=Span\{x^iy^j: 0 \leq i+j \leq n\}$
as a finite-dimensional representation; the l.h.s. of (A.12) annihilates
${\cal P}_n(x,y)$. If $y$ is a Grassmann variable,  $\tilde J^1_2 (n)$ is a
generator of the $q$-deformed superalgebra $osp(2,2)_q$ possessing
finite-dimensional representation.

Introducing two couples of the operators
$\tilde P_{1,2} , \tilde Q_{1,2}$ obeying the following relations:
\[ Q_{1}Q_{2}=qQ_{2}Q_{1}\ , \]
\[ P_{1} Q_{1}=1+q^2 Q_{1}P_{1}+(q^2-1) Q_{2}P_{2} \quad ,
\quad P_{1} Q_{2}=qQ_{2}P_{1} \ ,\]
\[ P_{2} Q_{1}=qQ_{1}P_{2}\quad , \quad P_{2} Q_{2}=1+q^2 Q_{2}P_{2} \ ,\]
\begin{equation}
 P_{1}P_{2} =q^{-1}P_{2}P_{1} \ .
\end{equation}
It is easy to derive a family of abstract
identities similar to (A.6) at $k=2$
\[ (Q_{1}^2 P_{1}\ +\ Q_{1}Q_{2} P_{2} - \{ n\} Q_{1})^{n+1} = \]
\begin{equation}
\sum_{k=0}^{n+1} q^{2n^2-n(k-2)+k(k-1)} {n+1 \choose k}_q
Q_{1}^{2n+2-k}Q_{2}^k
P_{1}^{n+1-k} P_{2}^k  \ ,
\end{equation}
while (A.12) is a particular case corresponding a special choice of a
representation (A.13) in  the form of (A.11).

Introducing a quantum hyperplane \cite{wz}, one can generalize
the whole family of operator identities (A.5) replacing first of all
continuous derivatives by finite-difference operators and then introducing
the abstract operators like it was done above for the case of (A.11)-(A.12).

The Lie-algebraic interpretation of above operator identities allows
to make a general conclusion that {\it an existence of finite-dimensional
representations of the Lie algebra $gl_k$ of differential operators
leads to an appearence of some specific operator identities analogous
to those families described above}.
\end{document}